\documentclass[12pt]{article}
\usepackage[tbtags]{amsmath}
\usepackage{amsfonts,amssymb,amsthm}
\usepackage{a4wide}
\usepackage{graphicx}
\usepackage{hyperref}
\usepackage{color}

\newtheorem{theorem}{Theorem}[section]
\newtheorem{lemma}[theorem]{Lemma}
\newtheorem{prop}[theorem]{Proposition}
\newtheorem{corollary}[theorem]{Corollary}

\theoremstyle{definition}
\newtheorem{definition}[theorem]{Definition}

\theoremstyle{remark}
\newtheorem{remark}[theorem]{Remark}

\begin{document}
\title{Rectangular Diagrams of Legendrian Graphs}
\author{Maxim Prasolov}
\date{}

\maketitle

\textsc{Abstract.}
In this paper Legendrian graphs in $(\mathbb{R}^3,\xi_{\mathrm{st}})$ are considered modulo Legendrian isotopy and edge contraction. To a Legendrian graph we associate a (generalized) rectangular diagram --- a purely combinatorial object. Moves of rectangular diagrams are introduced so that equivalence classes of Legendrian graphs and rectangular diagrams coincide. Using this result we prove that the classes of Legendrian graphs are in one-to-one correspondence with fence diagrams modulo fence moves introduced by Rudolph in \cite{Rud}, \cite{Rud2}.

\begin{center}
\textsc{Introduction.}
\end{center}

A Legendrian graph in $\mathbb{R}^3$ with standard contact structure $\xi_{\mathrm{st}}=\ker (dz+xdy)$ is a spatial graph  in $\mathbb{R}^3$ whose edges are Legendrian arcs, i.e. everywhere tangent to $\xi_{\mathrm{st}}$. In this paper Legendrian graphs are considered modulo continuous deformations in the class of Legendrian graphs (Legendrian isotopy) and modulo edge contraction. 

Legendrian graphs turn out to be useful in several problems. For example, they were used by Eliashberg and Fraser in \cite{EF} to classify topologically trivial Legendrian knots and in \cite{Giroux} by Giroux to establish a correspondence between contact structures and open book decompositions. In \cite{OP} O'Donnol and Pavelescu study Legendrian graphs as independent mathematical object.

In this paper we introduce a generalization of rectangular diagrams of links. In the literature oriented rectangular diagrams are commonly called "grid diagrams". Cromwell in \cite{Cromwell} and Dynnikov in \cite{Dyn} introduced rectangular elementary moves and proved that two rectangular diagrams represent isotopic links if and only if the diagrams are related by these elementary moves. Considering the equivalence relation generated only by some part of elementary moves one obtains a one-to-one correspondence between classes of rectangular diagrams and other geometric objects: 
\begin{itemize}
\item Legendrian links (\cite{OST});
\item transverse links (\cite{NT});
\item Birman-Menasco classes of braids (braids modulo conjugation and exchange move, \cite{NT}). 
\end{itemize}

The main result of this paper is a generalization of the first of these three correspondences. In \ref{mapG} we construct a map which associates to every rectangular diagram $R$ a Legendrian graph $G_R$. We prove the following:

\smallskip
\smallskip

\noindent\textbf{Theorem \ref{main_result}.} 
{\it The map $R\mapsto G_R$ induces a bijection between classes of generalized rectangular diagrams modulo elementary moves of type L and Legendrian graphs modulo Legendrian isotopy and edge contraction.
}
\smallskip
\smallskip

We apply this theorem to study fence diagrams modulo fence moves. They were introduced by Rudolph in \cite{Rud}, \cite{Rud2} to classify quasipositive surfaces. He asked: is it true that two fence diagrams are related by fence moves if and only if the corresponding quasipositive surfaces are isotopic? Baader and Ishikawa in \cite{BI} answered this question in the negative. Their idea was to construct a map from 3-valent Legendrian graphs modulo Legendrian isotopy to fence diagrams modulo fence moves. We prove the following:

\smallskip
\smallskip
\noindent\textbf{Corollary \ref{fence_leg}.} {\it A map introduced in \cite{BI} from 3-valent Legendrian graphs modulo Legendrian isotopy to fence diagrams modulo fence moves induces a bijection of Legendrian graphs modulo Legendrian isotopy and edge contraction with fence diagrams modulo fence moves.}
\smallskip
\smallskip

We prove this corollary by using a one-to-one correspondence between fence diagrams modulo fence moves and generalized rectangular diagrams modulo elementary moves of type L. We discuss this correspondence in section 4. So one can consider fence diagrams as a generalization of rectangular diagrams. It is also important to mention that rectangular diagrams were generalized for theta-graphs in \cite{DP}. The moves introduced in \cite{DP} are extented to the general case in this paper.

The organization of paper is as follows. In section 1 we recall known results on Legendrian graphs. In section 2 we introduce generalized rectangular diagrams and their moves. In section 3 we prove the main result. In section 5 we briefly concern a generalization of the correspondence between links and rectangular diagrams: we will show without details that two generalized rectangular diagrams are related by elementary moves if and only if the corresponding spatial graphs are related by isotopy and edge contraction.

\begin{center}
\textsc{Acknowledgements.}
\end{center}

I would like to thank my advisor Ivan Dynnikov for posing the problem. I am grateful to Masaharu Ishikawa and Sebastian Baader who shared with me a published version of their article. The author is partially supported by Laboratory of Quantum Topology of Chelyabinsk State University (Russian Federation government grant 14.Z50.31.0020)  and by science schools supporting grant NSh-4833.2014.1.

\section{Legendrian graphs}
\begin{definition}
By a {\it spatial graph} we mean a finite 1-dimensional CW-complex smoothly embedded in $\mathbb{R}^3$.
\end{definition}

In the definition of a spatial graph it is assumed that edges of the spatial graph are smooth simple arcs, and that at each vertex the tangent half-lines to incident edges are distinct.

\begin{definition}
A {\it Legendrian graph} is a spatial graph whose edges are tangent to the plane distribution defined by the kernel of the 1-form $$\alpha_{\mathrm{st}}=dz+xdy$$ and called {\it the standart contact structure}. We orient this distribution by the vector field $\frac{\partial}{\partial z}$ which is transverse to its planes. Two Legendrian graphs are called {\it Legendrian isotopic} if they are connected by a path in the space of Legendrian graphs.

Let us define a topology of the space of Legendrian graphs in detail.

Firstly, if we fix a 1-complex $G$ with finite number of vertices and edges and consider a space $\mathfrak{LE}(G)$ of Legendrian embeddings $G\hookrightarrow\mathbb{R}^3$ then this space possesses a compact-open topology.

Secondly, call two such Legendrian embeddings $\iota_k:G_k\hookrightarrow\mathbb{R}^3,\ k=1,2$ equivalent if there exists a {\it combinatorial equivalence} (a homeomorphism which maps bijectively vertices to vertices) $\varphi:G_1\stackrel{\approx}\to G_2$ such that the following diagram is commutative:
$$\begin{matrix}
G_1 & \stackrel{\iota_1}\hookrightarrow & \mathbb{R}^3 \\
\stackrel{\varphi}{\phantom{.}}\downarrow\phantom{.} & & \|\phantom{.} \\
G_2 & \stackrel{\iota_2}\hookrightarrow & \mathbb{R}^3
\end{matrix}$$

Thirdly, define a space of Legendrian graphs as $$\mathfrak{LG}=\bigsqcup\limits_{k=1}^{\infty}\mathfrak{LE}(G_k) / \sim,$$ where for each combinatorial class of graphs there is some representative $G_k$ and the relation $\sim$ is defined above. So two Legendrian graphs are Legendrian isotopic iff they lie in the same component of $\mathfrak{LG}$.

Note that two Legendrian isotopic Legendrian graphs are not neccesary ambient isotopic. The reason is that a diffeomorphism preserves linear relations on tangent vectors to edges at each vertex.

\end{definition}

Note that at each vertex of a Legendrian graph the tangent spaces to incident edges lie in a oriented plane of the distribution. This gives a cyclic order of edges incident to this vertex. So every Legendrian graph is a ribbon graph.

\begin{figure}
\begin{center}
\includegraphics{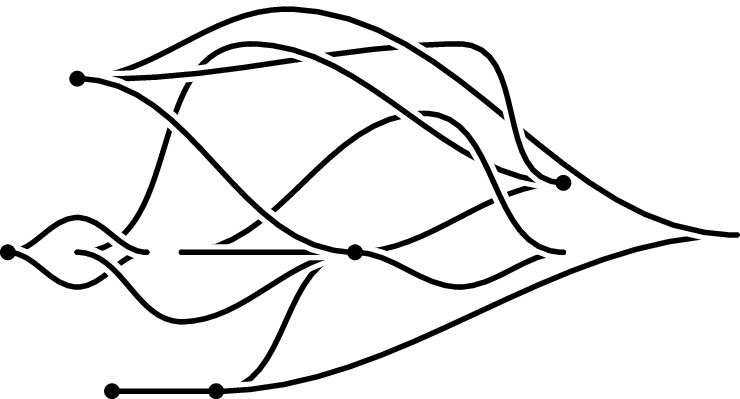}
\end{center}
\caption{Front of a Legendrian graph. Crossing information is redundant since an undercrossing branch always has greater slope.}
\label{front_example}
\end{figure}

It is convenient to specify Legendrian graphs by drawing their projections to $yz$-plane, which are called {\it front projections} or simply {\it fronts}. The front is a collection of piecewise-smooth curves whose singularities have the form of a cusp and whose tangents are never vertical (see Fig. \ref{front_example}). At self-intersections and cusp-points (for visual evidence only) we display which branch is overcrossing, and which is undercrossing. A generic front projection uniquely determines the corresponding Legendrian curve in $\mathbb{R}^3$ since $x$-coordinate of any point of the curve can be recovered from the relation $x=-dz/dy.$ So the information about overcrossing and undercrossing is unneccesary.

\begin{figure}
\center{\includegraphics{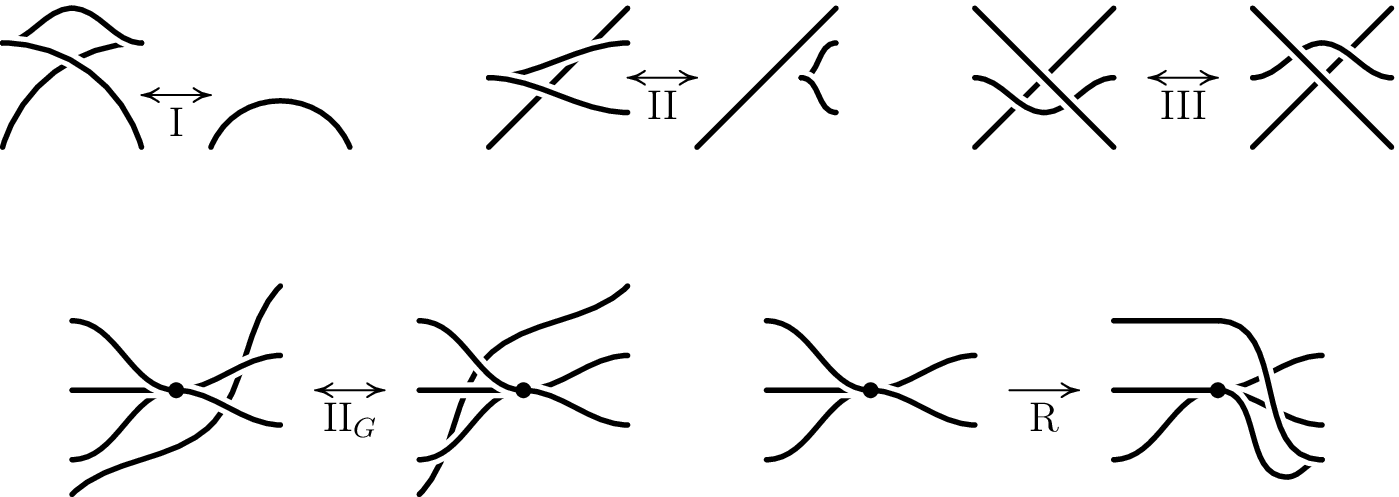}}
\caption{Moves of fronts. We should also add moves obtained by reflection with respect to horizontal or vertical axis (with overcrossing and undercrossing exchanged) and by $\pi$-rotation. In the moves $\mathrm{II}_G$ and R number of edges pointing to the left and to the right may be arbitrary. }
\label{front_moves}
\end{figure}

\begin{theorem}[\cite{BI}]\label{BItheo}
Two generic fronts represent Legendrian isotopic Legendrian graphs iff they are related by moves which are illustrated in Figure \ref{front_moves}.
\end{theorem}

Moves I,II and III generate Legendrian isotopy equivalence between Legendrian links (see \cite{Swia} for the proof), and the rest moves engage the vertices of the graph. 

\begin{figure}
\center{\includegraphics{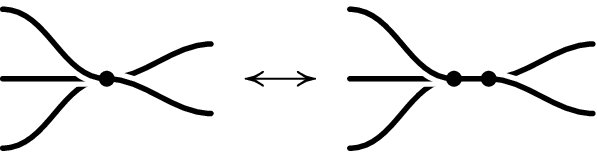}}
\caption{Blow-up and edge contraction}
\label{edge_contraction}
\end{figure}

\begin{definition}
Consider one more front move: an {\it edge contraction} (see Fig.~\ref{edge_contraction}). Suppose we have two vertices connected by a short horizontal segment on the front such that the rest edges which are incident to the left (respectively, right) vertex emerge to the left (respectively, right). An edge contraction consists in contracting the small edge to a vertex such that the order and the direction of emergence of the remaining edges preserve. The inverse move is called {\it blow-up}.
\end{definition}

The following theorem is the unique such result of this section that will be used in the following sections.

\begin{theorem}
\label{front_moves_reduced}
Legendrian graphs modulo Legendrian isotopy and edge contraction are in one-to-one correspondence with generic fronts modulo moves R, $\mathit{II}_G$, III and edge contraction.
\end{theorem}

\begin{figure}
\center{\includegraphics{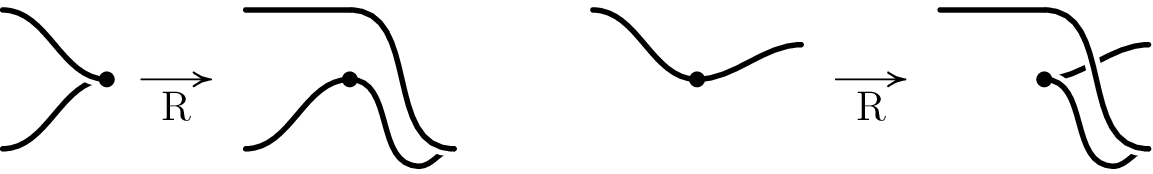}}
\caption{On the left: moving 2-valent vertex through the cusp using R-move. On the right: applying move I using R-move.}
\label{r_move_application}
\end{figure}

\proof
To prove this theorem we represent the moves I and II as a combination of other moves. This suffices by theorem~\ref{BItheo}.

Consider the move I. Edge considered in the move is connected to some vertex. Make a blow-up at this vertex and move new (2-valent) vertex (using several moves $\mathrm{II}_G$ to overpass crossings and R-moves to overpass cusps, see Fig. \ref{r_move_application}) to the place where we want to apply move I. Then apply R-move at this vertex, see Fig. \ref{r_move_application}. Then we want to remove this vertex. Draw it back (again using moves $\mathrm{II}_G$ and R) to adjacent vertex and apply edge contraction. So move I is a combination of moves R, $\mathrm{II}_G$, edge contraction and blow-up. Similarly move II is a combination of other moves. 
\endproof

In the end of this section we introduce a slightly generalized version of the blow-up and edge contraction moves which will be useful in the following sections, and we introduce an example of move that is not a blow-up.

We introduced a blow-up and an edge contraction as moves of fronts only. But in fact these operations are geometric. Choose some vertex and denote edges, which are incident to this vertex, say, by $a_1,\dots,a_n,b_1,\dots,b_m$ in the cyclic order. So we just divided edges into two groups $\{a_i\}_{i=1}^n$ and $\{b_j\}_{j=1}^m$ which are not "interleaving" in the cyclic order. To make a (geometric) blow-up one should substitute for the vertex a small (Legendrian) edge $e$ and slightly move edges $a_1,...,a_n,b_1,...b_m$ such that at the first end of $e$ edges are arranged in the cyclic order as $a_1,\dots,a_n,e$ and in the second end as $e,b_1,\dots,b_m.$

\begin{figure}
\center{\includegraphics{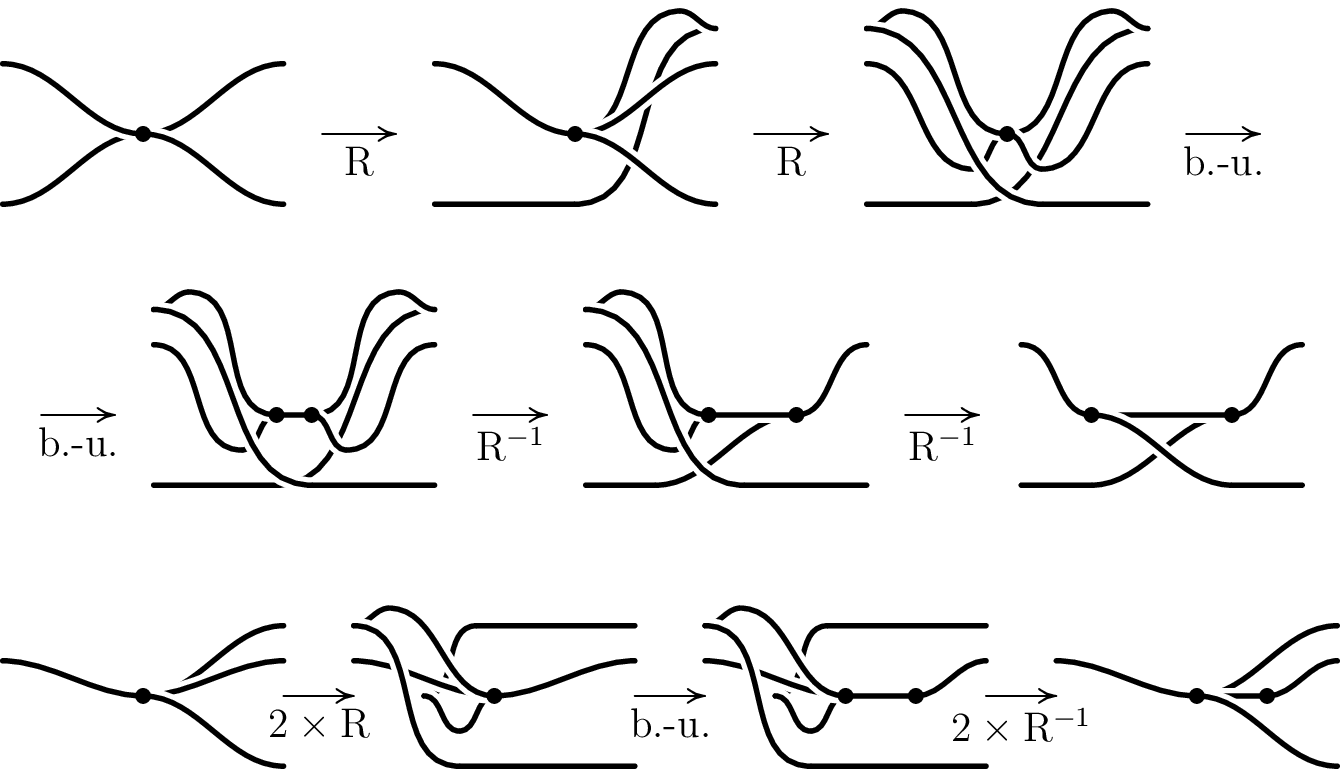}}
\caption{Blow-up variations}
\label{generalized_blow-up}
\end{figure}

By combining move R, its inverse and blow-up front move one can perform a geometric blow-up. In the Figure \ref{generalized_blow-up} each of four branches emerging from the vertex can be substituted by any number (even zero) of branches. In this Figure all cases of geometric blow-ups are considered up to rotation by $\pi$. 

\begin{figure}
\center{\includegraphics{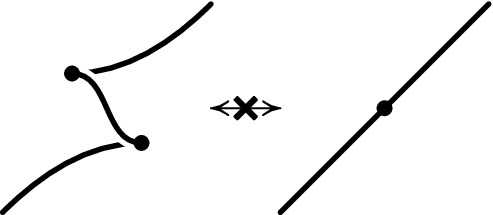}}
\caption{Wrong interpretation of edge contraction}
\label{forbidden_move}
\end{figure}

One can be confused that the move in Figure \ref{forbidden_move} may be regarded as a blow-up or an edge contraction. In fact, in some cases this move leads to a graph which can not be obtained by Legendrian isotopy, blow-ups and edge contractions from the initial graph. To demonstrate this we introduce one more definition.

\begin{definition}
\label{ribbon}
Fix some Legendrian graph and consider a surface (with boundary) embedded in $\mathbb{R}^3$ which contains in its interior this graph and is tangent to the contact structure at all points of the graph. If this surface is sufficiently small then its isotopy class is well defined. Also this class is preserved under Legendrian isotopy of the graph and edge contraction or blow-up. Call this small surface a {\it ribbon} of the Legendrian graph.
\end{definition}

In the case of the graph with two vertices connected by two edges the move in Figure \ref{forbidden_move} alters the linking number of the ribbon boundary components  which is an annulus. By the way, in this case the graph represents a knot and this linking number is the Thurston-Bennequin invariant of the Legendrian knot.

\section{Generalized rectangular diagrams}

\begin{figure}
\center{\includegraphics{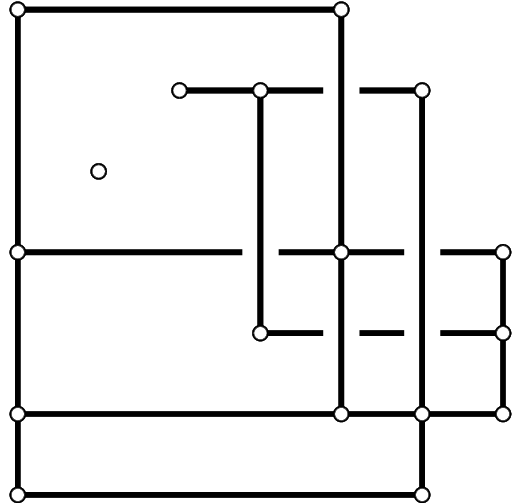}}
\caption{Example of generalized rectangular diagram}
\label{diagram}
\end{figure}

\begin{definition}
\label{GRD}
A {\it generalized rectangular diagram} is some finite subset of the plane. Elements of this set are called {\it vertices} of the rectangular diagram. If some vertices lie on the same horizontal or vertical line then we connect these vertices by a segment with ends in extreme points. We call this segment an {\it edge} of the rectangular diagram. 

We will assume that two diagrams are the same if they are {\it combinatorially equivalent}: if there exists two increasing functions $f,g:\mathbb{R}\mapsto\mathbb{R}$ such that one diagram maps to another by a map $(x,y)\mapsto(f(x),g(y))$.

If intersection of two edges is not a vertex then the vertical line is overcrossing and the horizontal line --- undercrossing. Such diagram should be interpreted as a planar diagram of a spatial graph.
\end{definition}

An example of generalized rectangular diagram is shown in Figure \ref{diagram}. We allow edges which contain one vertex.

Now we will introduce {\it elementary moves} of generalized rectangular diagram. We define each move by specifying positions of new vertices only. We do not specify positions of edges because they are uniquely determined by positions of vertices (edges are the segments which connect extreme vertices on the horizontal and vertical lines).

\begin{figure}[h]
\center{\includegraphics{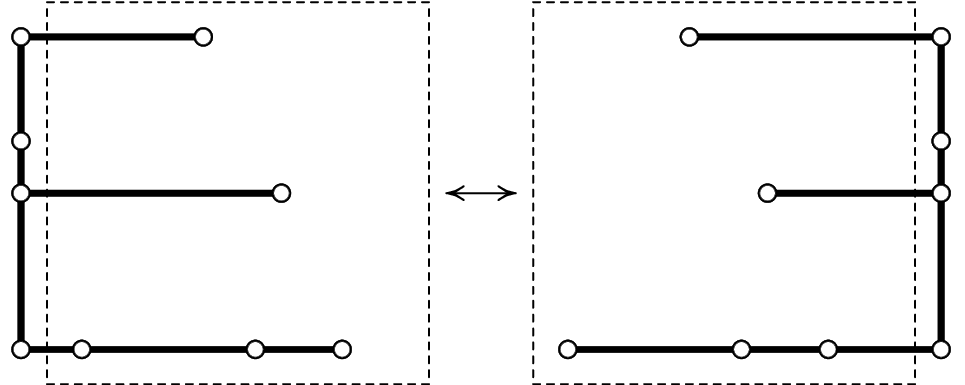}}
\caption{Cyclic permutation}
\label{cyclic_permutation}
\end{figure}

\begin{itemize}
\item A {\it cyclic permutation} of vertical or horizontal edges consists in moving one of the extreme (top, bottom, left, or right) edges onto the opposite side, see Fig. \ref{cyclic_permutation}. Only one of two coordinates of vertices on the edge alters.

\begin{figure}
\center{\includegraphics{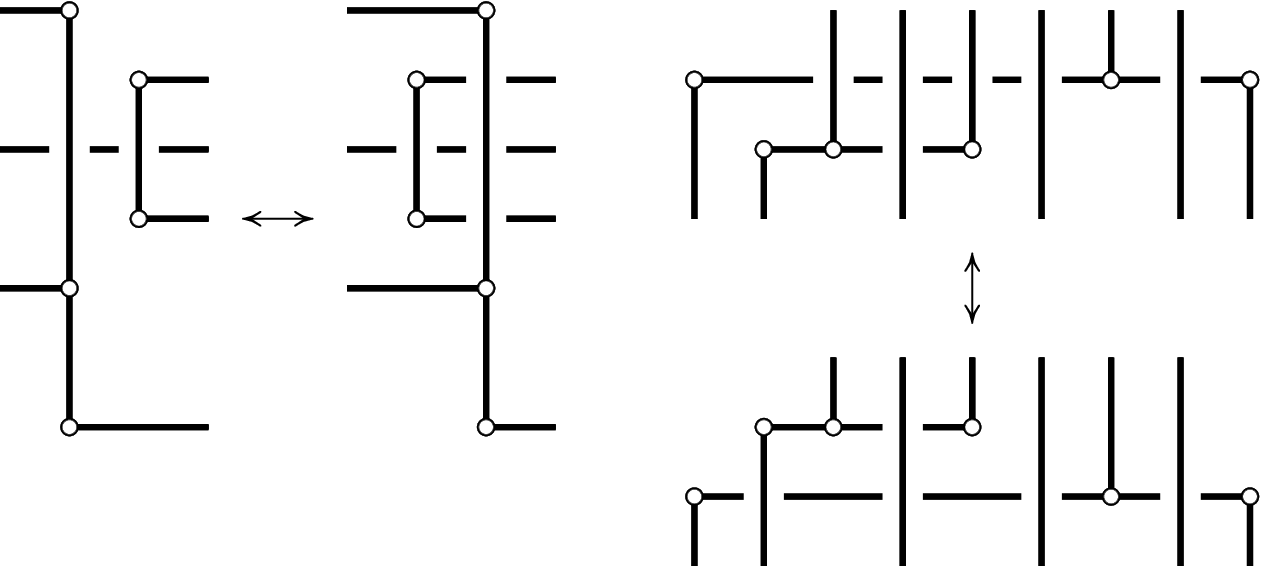}}
\caption{Commutations}
\label{commutation}
\end{figure}

\item A vertical (respectively, horizontal) {\it commutation} consists in exchange of the horizontal (respectively, vertical) positions of two neighboring vertical (respectively, horizontal) edges provided that these edges do not {\it interleave}, i.e. in the cyclic order the $y-$(respectively, $x-$)coordinate of vertices on one edge is smaller than of vertices on the other edge. The edges are regarded neighboring if there are no vertices of the diagram between the parallel straight lines containing these edges, see Figure \ref{commutation}.

\begin{figure}
\center{\includegraphics{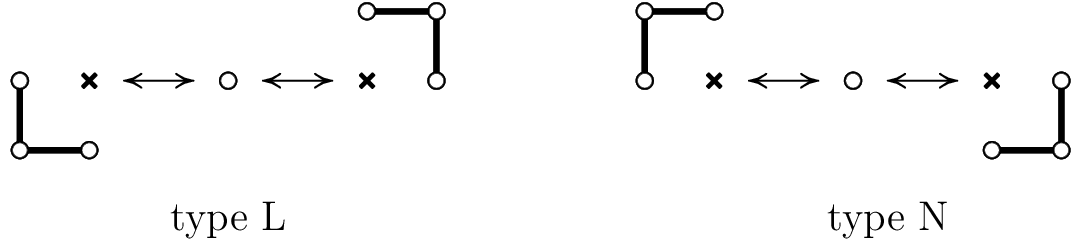}}
\caption{(De)stabilizations of two types}
\label{stabilizations}
\end{figure}

\item A {\it stabilization} consists in replacement of a vertex by three new ones that together with the deleted one form the vertex set of a small square and addition of two short edges that are sides of the square; the inverse operation is called a {\it destabilization} (see Fig. \ref{stabilizations}). 
If two new edges emerge from the stabilization point from their common end downward and leftward or upward and rightward, then the stabilization is of type L (i.e. "Legendrian"), and otherwise of type N (i.e. "non-Legendrian").

\begin{figure}
\center{\includegraphics{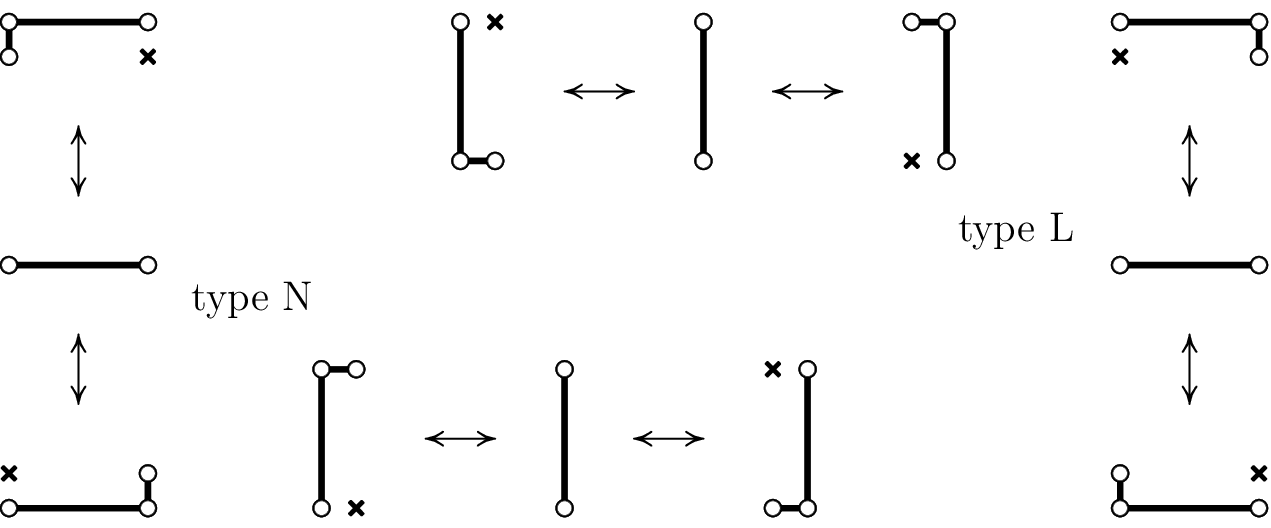}}
\caption{End shift}
\label{end_shift}
\end{figure}

\item Let $v_1$ and $v_2$ be adjacent vertices on some edge and $v_2$ has greater (respectively, smaller) coordinate than $v_1$. An {\it end shift of type L (respectively, N)} consists in introducing two vertices which are shifts of $v_1$ and $v_2$ by $\varepsilon>0$ in the positive direction perpendicular to the edge and removing $v_1$ (see Fig. \ref{end_shift}). Or it consists in introducing two vertices which are shifts of $v_1$ and $v_2$ by $\varepsilon$ in the negative direction perpendicular to the edge and removing $v_2$ (see Fig. \ref{end_shift}). Number $\varepsilon$ is chosen sufficiently small such that there are no vertices between line connecting $v_1$ with $v_2$ and line connecting shifts of these vertices. The inverse operation to end shift (of type L) is a combination of (type L) end shift, commutations and (type L) destabilization. See Fig. \ref{end_shift_inverse}

\begin{figure}[h]
\center{\includegraphics{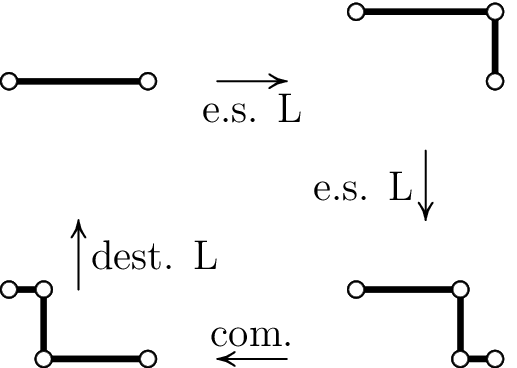}}
\caption{Inversing an end shift}
\label{end_shift_inverse}
\end{figure}

\item {\it A vertex addition} consists in introducing a new vertex on some vertical (respectively, horizontal) line containing some another vertex provided that there are no vertices on the same horizontal (respectively, vertical) line containing this new vertex.
\end{itemize}

A cyclic permutation, a commutation, a (de)stabilization and an end shift are straightforward generalizations of the moves with the same name introduced in \cite{DP} for theta-graphs.

For convenience, by {\it elementary moves of type L (respectively, N)} we will mean a cyclic permutation, commutations, (de)stabilizations of type L (respectively, N), end shift of type L (respectively, N) and a vertex addition.

\begin{theorem}
\label{elementary_basis}
Equivalence relation on generalized rectangular diagrams generated by elementary moves of type L (respectively, N) in fact is generated by commutations, an end shift of type L (respectively, N) and a vertex addition.
\end{theorem}
\proof
We show that a stabilization and a cyclic permutation are combinations of other moves.

\begin{figure}[h]
\center{\includegraphics{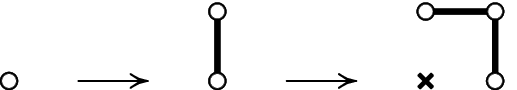}}
\caption{Stabilization as a combination of one vertex addition and one end shift}
\label{stabilization_combine}
\end{figure}

Stabilization of type L (respectively, N) is a combination of a vertex addition and a type L (respectively, N) end shift: add a vertex near stabilization point and apply an end shift to the pair of the initial vertex and the new one (see Fig. \ref{stabilization_combine}).

\begin{figure}[h]
\center{\includegraphics{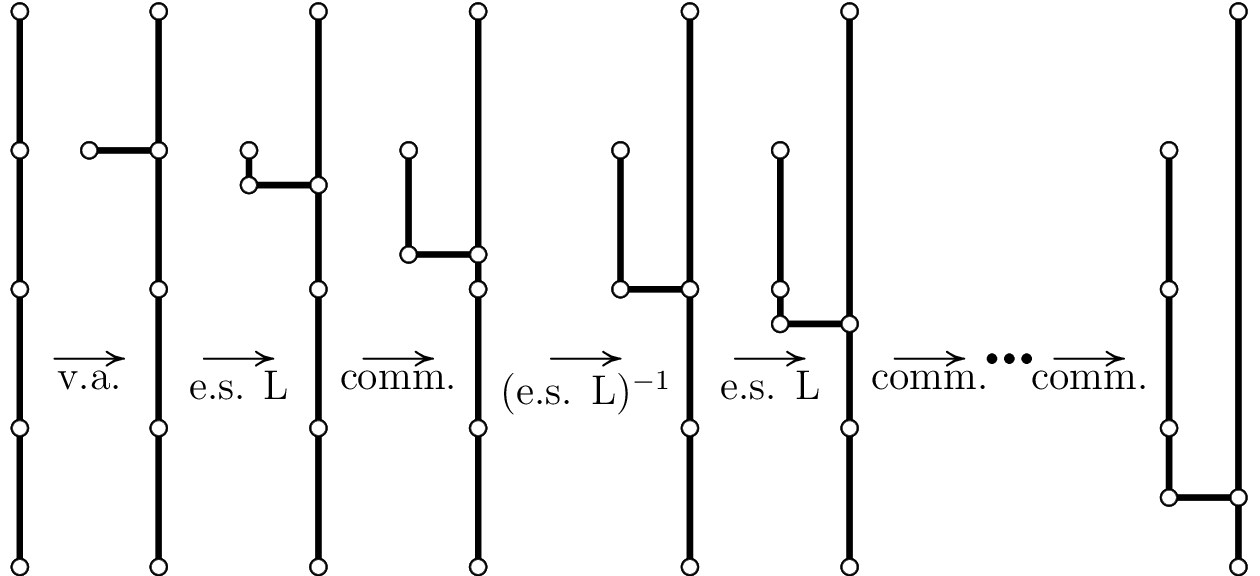}}
\caption{Edge breaking is a combination of elementary moves}
\label{edge_breaking_combine}
\end{figure}

{To consider the case of the cyclic permutation we introduce an auxiliary move called edge breaking.} Horizontal (respectively, vertical) {\it edge breaking of type L} consists in separating vertices of a horizontal (respectively, vertical) edge into two parts forming neighboring commuting edges $e_1$ and $e_2$, suppose $e_2$ has bigger vertical (respectively, horizontal) coordinate than $e_1$ has, and jointing them by a short edge in such a way that in a cyclic order a horizontal (respectively, vertical) coordinate of short edge is smaller than of vertices of $e_1$ but bigger than of $e_2$. To specify an edge breaking of type N switch $e_1$ and $e_2$ in the above definition.
\begin{figure}[h]
\center{\includegraphics{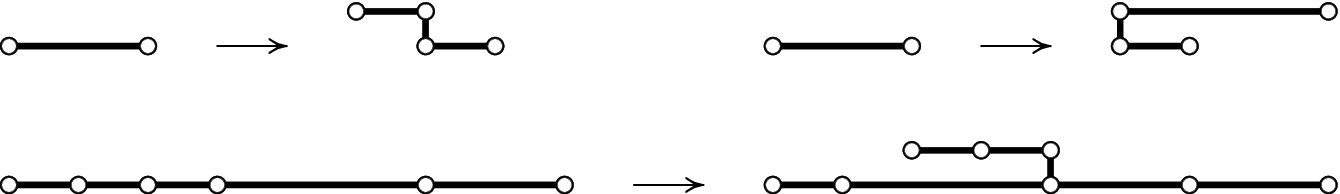}}
\caption{Examples of horizontal type L edge breaking. To obtain examples of vertical type L edge breaking reflect the picture with respect to diagonal line. To obtain examples of type N edge breaking reflect the picture with respect to horizontal or vertical lines.}
\label{edge_breaking}
\end{figure}
An example is shown in Fig.\ref{edge_breaking}. Edge breaking of type L (respectively, N) is a combination of type L (respectively, N) elementary moves. Proof is shown in Fig. \ref{edge_breaking_combine} for type L: recall that the inverse of a type L end shift, used in Figure, is a combination of type L elementary moves. Reflecting this picture with respect to some horizontal line gives the proof for the type N. The inverse operation is called {\it edge jointing}. This operation is also a combination of elementary moves.

\begin{figure}[p]
\center{\includegraphics{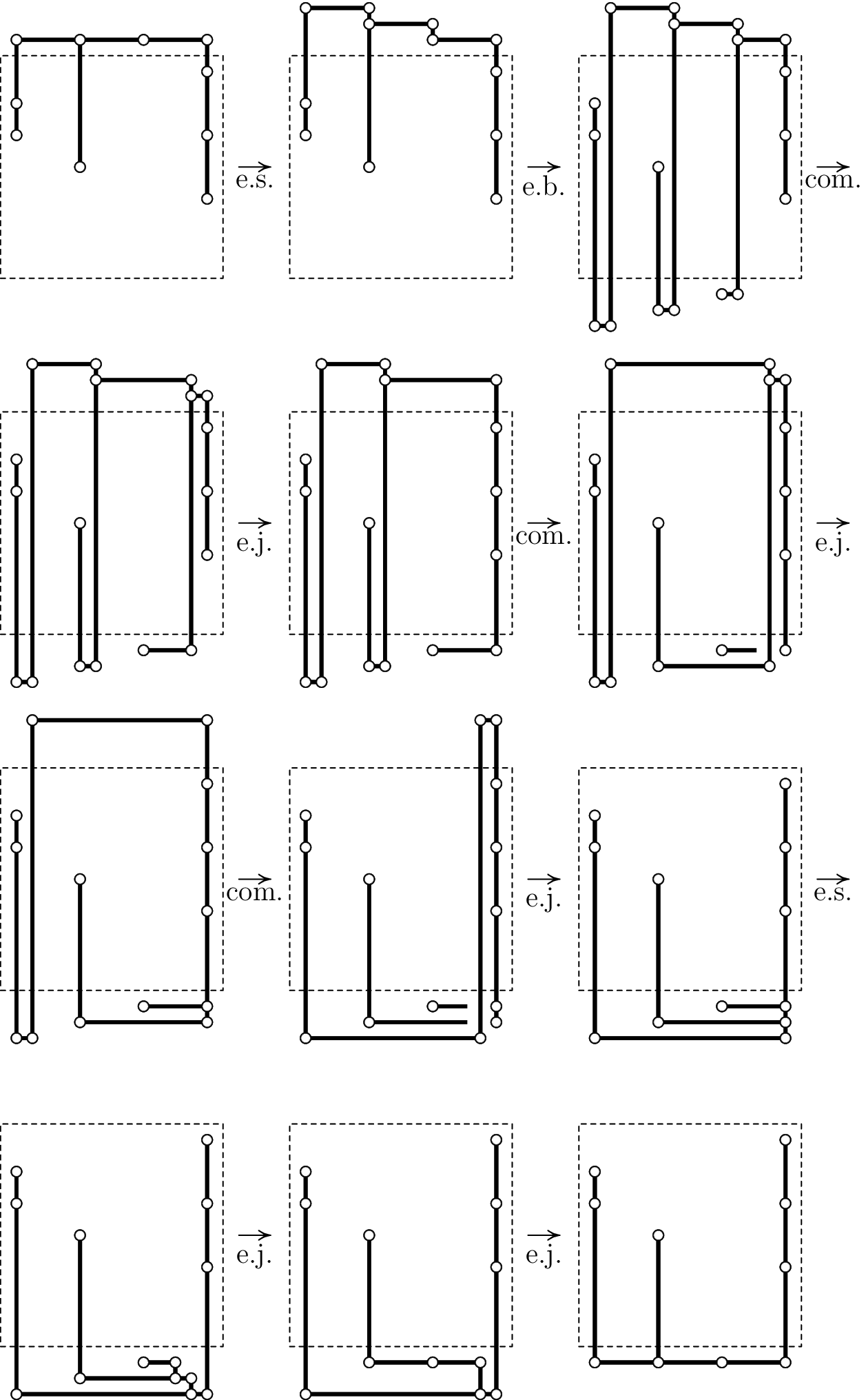}}
\caption{Cyclic permutation is a combination of other (type L) elementary moves}
\label{cyclic_permutation_combine}
\end{figure}

Now we can see in Fig.~\ref{cyclic_permutation_combine} that a cyclic permutation is a combination of other type L elementary moves.

\endproof

\begin{remark} A vertex addition is a combination of a stabilization, an end shift and a destabilization (see Fig.~\ref{add_vertex_combine}).
\begin{figure}[h]
\center{\includegraphics{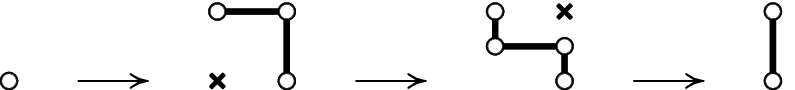}}
\caption{A vertex addition as a combination of (de)stabilizations and one end shift}
\label{add_vertex_combine}
\end{figure}
\end{remark}

\subsection{Flype}
In this subsection we extend to generalized rectangular diagrams a {\it flype} --- a move introduced by Dynnikov in \cite{Dyn2} for rectangular diagrams. We do not use this move in the following sections, so if you are interested only in the main result of this paper you can freely skip this subsection.

\begin{figure}
\center{\includegraphics{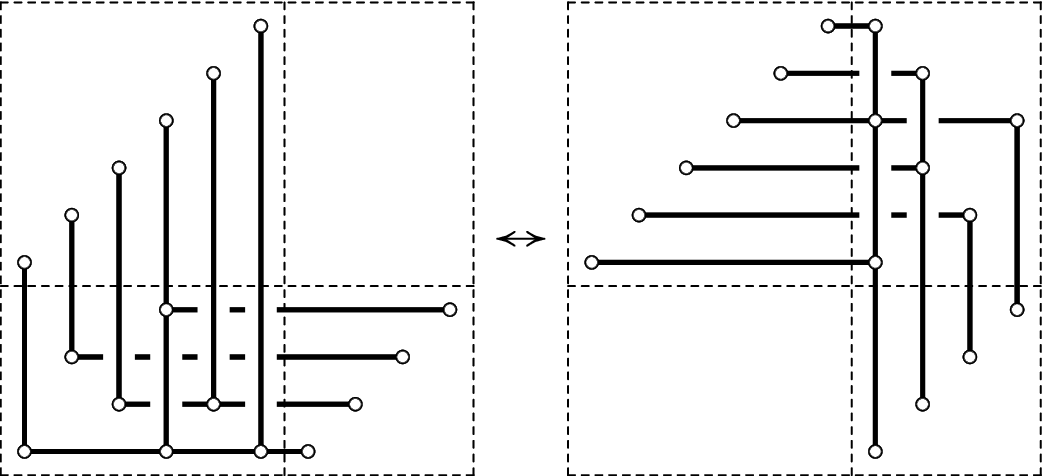}}
\caption{An example of (type L) flype}
\label{flype}
\end{figure}

\begin{definition}
Let $x_0<x_1<x_2$, $y_0<y_1<y_2$ and $$r_{ij}:=\{(x,y) \mid x_i<x<x_{i+1}, y_j<y<y_{j+1}\} \text{ for } i,j=0,1.$$ Let $R$ be a (generalized) rectangular diagram such that:
\begin{itemize}
\item the rectangle $r_{11}$ does not contain any vertices of $R$;
\item for each rectangle $r_{01}$ and $r_{10}$ for any pair of vertices lying in the rectangle the right vertex is higher than the left one;
\item for every vertex $v$ of $R$ lying in the rectangle $r_{00}$ there is a pair of vertices $v_{01}$ and $v_{10}$, such that $v_{01}$ has the same $x$-coordinate as $v$ and lies in $r_{01}$, and $v_{10}$ has the same $y$-coordinate and lies in $r_{10}$.
\end{itemize}
A {\it flype} of type L consists in moving each vertex $v$ contained in $r_{00}$ to the fourth corner of a rectangle with corners $v$, $v_{01}$ and $v_{10}$.

The inverse move is also called flype of type L. A move obtained by a conjugation by a horizontal symmetry is called {\it flype of type N}.

An example is shown in Fig.~\ref{flype}.
\end{definition}

\begin{prop}
\label{flype_combination}
A flype of type L is a combination of type L elementary moves.
\end{prop}
\proof
Order vertices of $R$ lying in $r_{00}$ from top to bottom and from right to left. Note that a vertex $v$ in $r_{00}$ can be moved to the fourth corner of a rectangle with corners $v, v_{01},v_{10}$ by one end shift, commutations and one inverse of an end shift. So apply these moves for each vertex in turn beginning from the smallest vertex.
\endproof

\begin{figure}
\center{\includegraphics{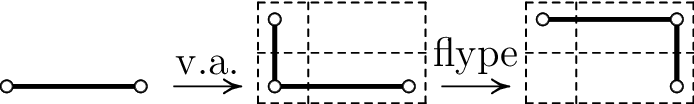}}
\caption{An end shift is a combination of a flype and a vertex addition}
\label{end_shift_flype}
\end{figure}

\begin{prop}
An end shift of type L is a combination of a vertex addition and a type L flype. A stabilization of type L is a combination of two vertex additions and a type L flype. A commutation is a combination of two type L flypes and several additions and removings of vertices.
\end{prop}
\proof Cases of an end shift and a stabilization are similar and very easy. See Fig.~\ref{end_shift_flype}.

Now to the case of a commutation. Suppose our commuting edges are vertical. The horizontal case can be obtained by a reflection in a line $x=y.$

\begin{figure}
\center{\includegraphics{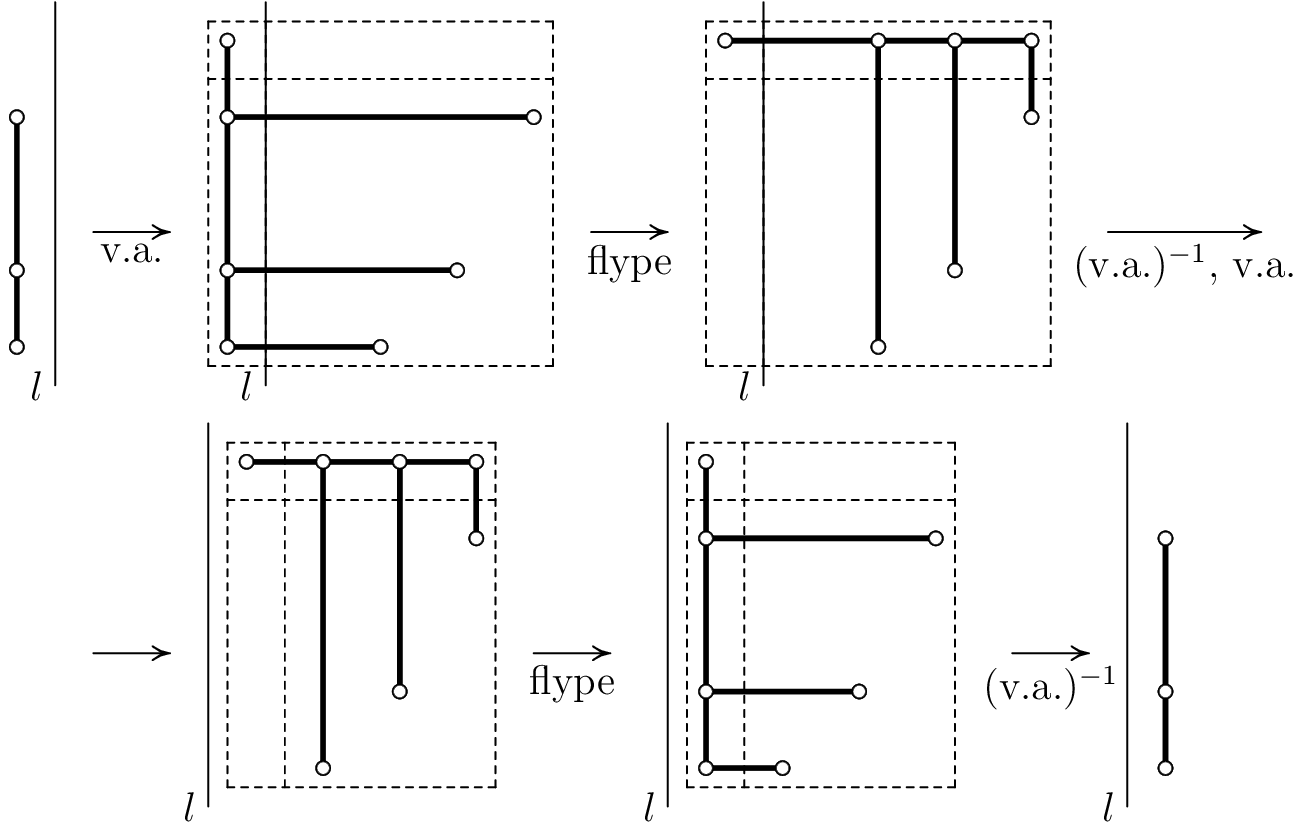}}
\caption{A commutation is a combination of two flypes, vertex additions and removings}
\label{commutation_flype}
\end{figure}

There are two types of a commutation: when projections to a $y-$axis of edges do not intersect and when a projection of the first edge contains a projection of the second edge. In the latter case choose the second edge and in the first case choose any edge. Denote the chosen edge by $e$. Denote a line containing the other edge by $l$.

Is is sufficient to consider a case when the edge $e$ lies to the left from the line $l$, since a commutation in the other case is the inverse move. Choose a rectangle $r$ such that:
\begin{itemize}
\item its sides are parallel to the $x-$ and $y-$ axes;
\item the line $l$ intersects an interior of $r$;
\item $r$ contains $e$ and any vertex contained in $r$ lies on $e$.
\end{itemize}

Let $l'$ be a horizontal line which is above $e$ but intersects an interior of $r$. Then $l$ and $l'$ divide $r$ into four pieces $r_{00}, r_{01}, r_{10}, r_{11}$ ($r_{00},r_{01}$ are two left pieces and $r_{00}, r_{10}$ are two bottom pieces). Add a vertex in $r_{10}$ for each vertex of $e$ and add a vertex $v$ in $r_{01}$ on the horizontal level of $e$ such that it is possible to perform a flype. Perform the flype and remove $v$: now all vertices of $e$ are moved to $r_{11}$. 

Introduce a vertical line $l''$ which intersects an interior of $r_{11}$ and such that there are no vertices in $r_{11}$ or $r_{10}$ to the left of $l''$. The line $l''$ divides $r_{10}$ into two pieces. Denote them by $r'_{00}$ for the left one and by $r'_{10}$ for the right one. Similarly $r_{11}=r'_{01}\cup r'_{11}$. Add a vertex in $r'_{01}$ to make it possible to perform a flype with respect to the rectangles $r'_{ij}$, $i,j=0,1$. After that remove all vertices in $r'_{01}\cup r'_{10}$ and we are done. See Fig.\ref{commutation_flype}.
\endproof

\begin{corollary}
The equivalence relation on generalized rectangular diagrams modulo elementary moves of type L is generated by flypes of type L and a vertex addition.
\end{corollary}
\proof This follows from the previous proposition and Theorem~\ref{elementary_basis}.
\endproof

\begin{prop}
The reflection in the line $x=y$ of a rectangular diagram is a combination of type L elementary moves.
\end{prop}
\proof
Choose rectangles $r_{00},r_{01},r_{10},r_{11}$ as in the definition of a flype such that $r_{00}$ contains the whole rectangular diagram. For each horizontal edge of the diagram add a vertex in $r_{10}$ and for each vertical edge of the diagram add on the line containing the edge a vertex in $r_{01}$ such that it is possible to make a flype. After making a flype remove all vertices in $r_{01}\cup r_{10}$. The obtained diagram is combinatorially equivalent to the diagram obtained by the reflection in the line $x=y$ from the initial diagram. By Proposition~\ref{flype_combination} --- that a type L flype is a combination of type L elementary moves --- we are done.
\endproof

\newpage
\section{Correspondence theorem}
\subsection{Correspondence map}

\begin{figure}
\center{\includegraphics{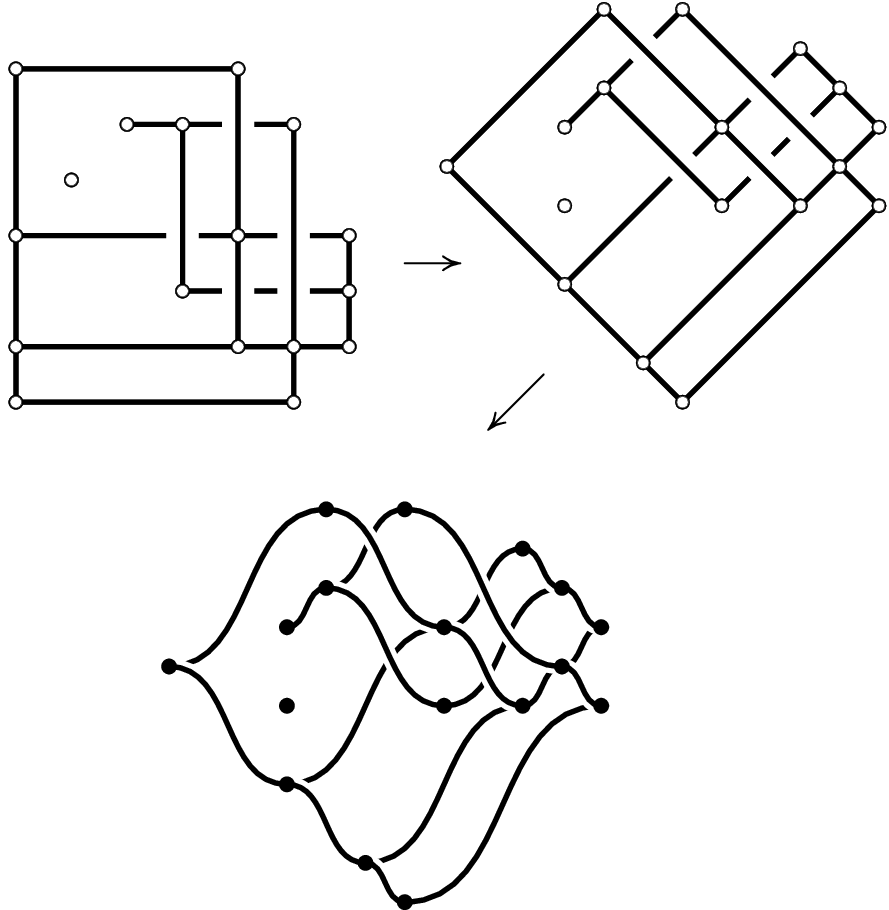}}
\caption{Front of a Legendrian graph obtained from generalized rectangular diagram}
\label{rect_front}
\end{figure}
\begin{definition}
\label{mapG}
Let $R$ be a generalized rectangular diagram. Rotate this diagram by $\pi/4$ counterclockwise, smooth out the corners pointing up and down and turn into cusps corners pointing to the left or to the right. Turn all vertices of the diagram to vertices of a graph specified by the the obtained front. Denote this graph by $G_R$. See Figure \ref{rect_front} for an example. 
\end{definition}

\begin{theorem}
\label{main_result}
The map $R\mapsto G_R$ induces a bijection between classes of generalized rectangular diagrams modulo elementary moves of type L and Legendrian graphs modulo Legendrian isotopy and edge contraction / blow-up.
\end{theorem}

A proof of this theorem decomposes into two parts.

The first part is to check that the correspondence map is well defined and the second --- to construct the inverse map.

Now to the first part of the proof. So we want to demonstrate that if $R'$ is obtained from $R$ by an elementary move of type L then $G_{R'}$ and $G_{R}$ are equivalent modulo Legendrian isotopy and edge contraction / blow-up.
By Theorem \ref{elementary_basis} --- that elementary moves are generated by commutations, end shifts and vertex additions --- it is sufficient to consider the case when $R'$ is obtained from $R$ by a commutation, an end shift of type L or a vertex addition.

\begin{figure}
\center{\includegraphics{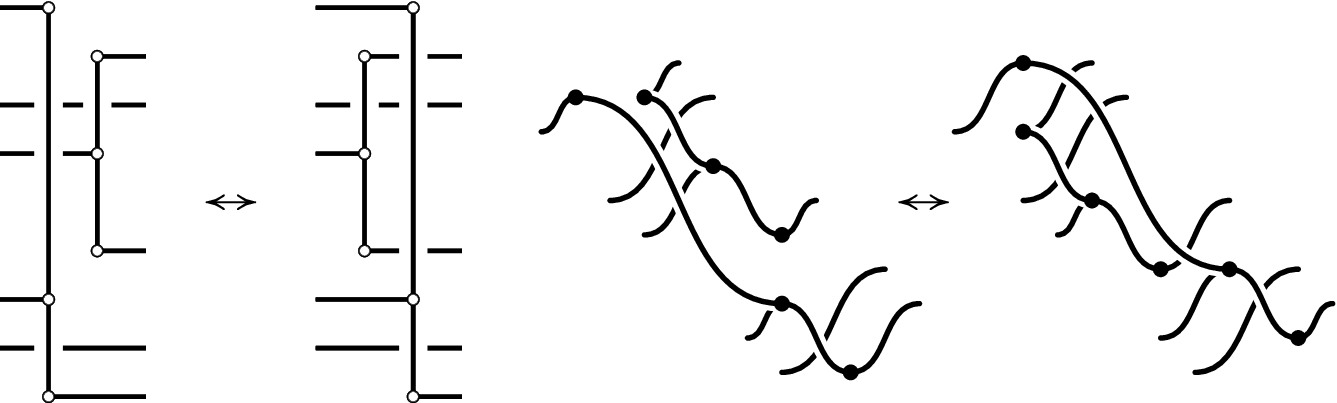}}
\caption{Equivalence of fronts in the case of a commutation}
\label{equiv_front_commutation}
\end{figure}
\noindent\textbf{Commutation.} If projections of two commuting edges to parallel (to these edges) axis do not intersect then fronts of $G_{R'}$ and $G_R$ are isotopic. Otherwise one edge is strictly smaller than another one. Each vertex lying on the smaller edge leads to move $\mathrm{II}_G$ of the fronts and each crossing on this edge leads to move $\mathrm{III}$. See Fig. \ref{equiv_front_commutation}.

\begin{figure}
\center{\includegraphics{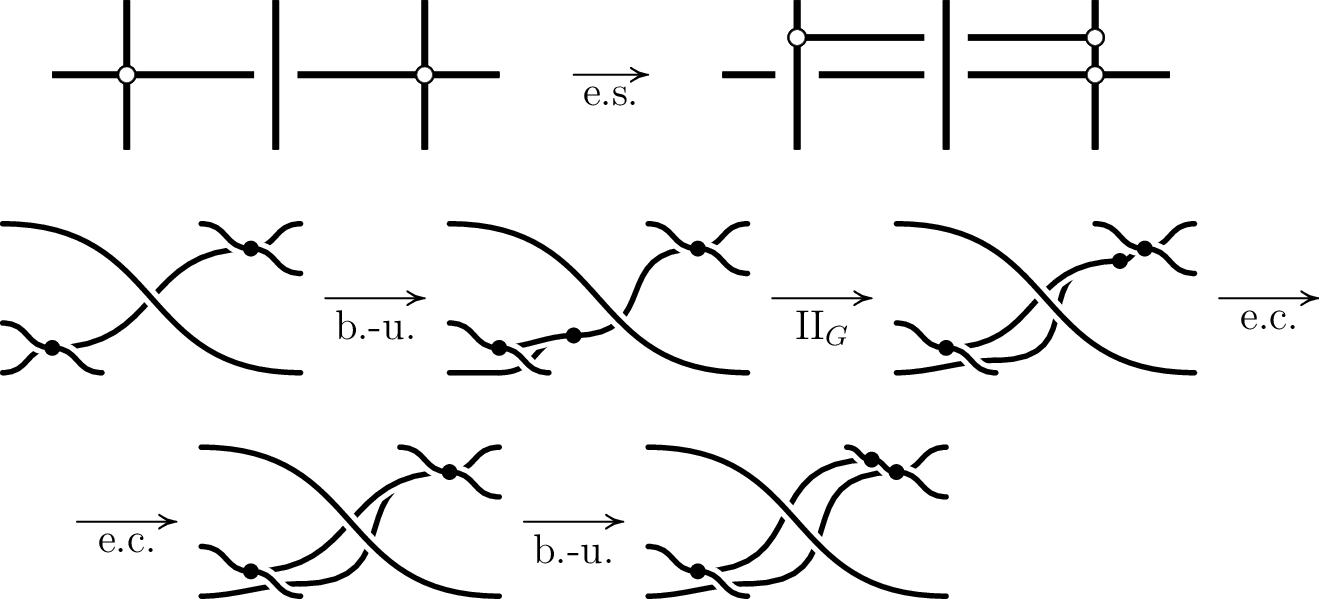}}
\caption{Equivalence of fronts in the case of an end shift of type L}
\label{equiv_front_end_shift}
\end{figure}
\noindent\textbf{End shift of type L.} In this case one front can be obtained from the other by a combination of blow-ups, edge contractions and moves $\mathrm{II}_G$ (see Fig. \ref{equiv_front_end_shift}). 

\noindent\textbf{Vertex addition.} In this case it is evident that the graphs are related by isotopy and a blow-up.

\subsection{Construction of the inverse map}
\begin{definition}[\cite{BI}]
\label{backslash_position}
A front projection is called {\it in backslash position} if all its tangent lines lie in $(\pi/2,\pi)\cup(3\pi/2,2\pi).$
\end{definition}

\begin{figure}[h]
\center{\includegraphics{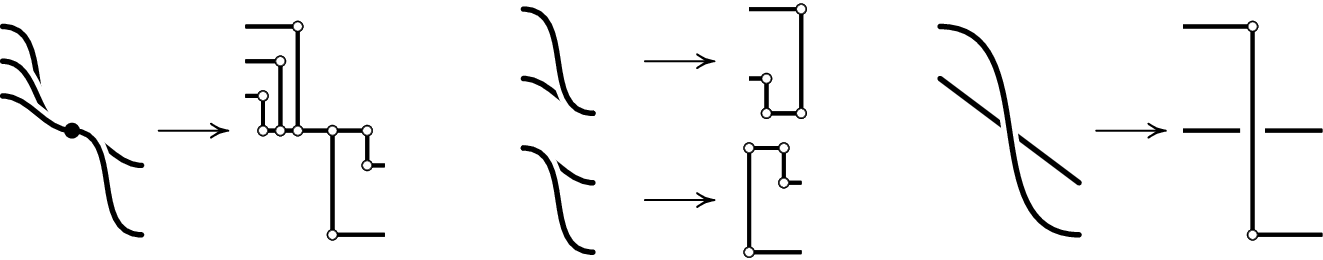}}
\caption{Parts of a generalized rectangular diagram approximating a front in backslash position}
\label{approx_parts}
\end{figure}
\begin{figure}
\center{\includegraphics{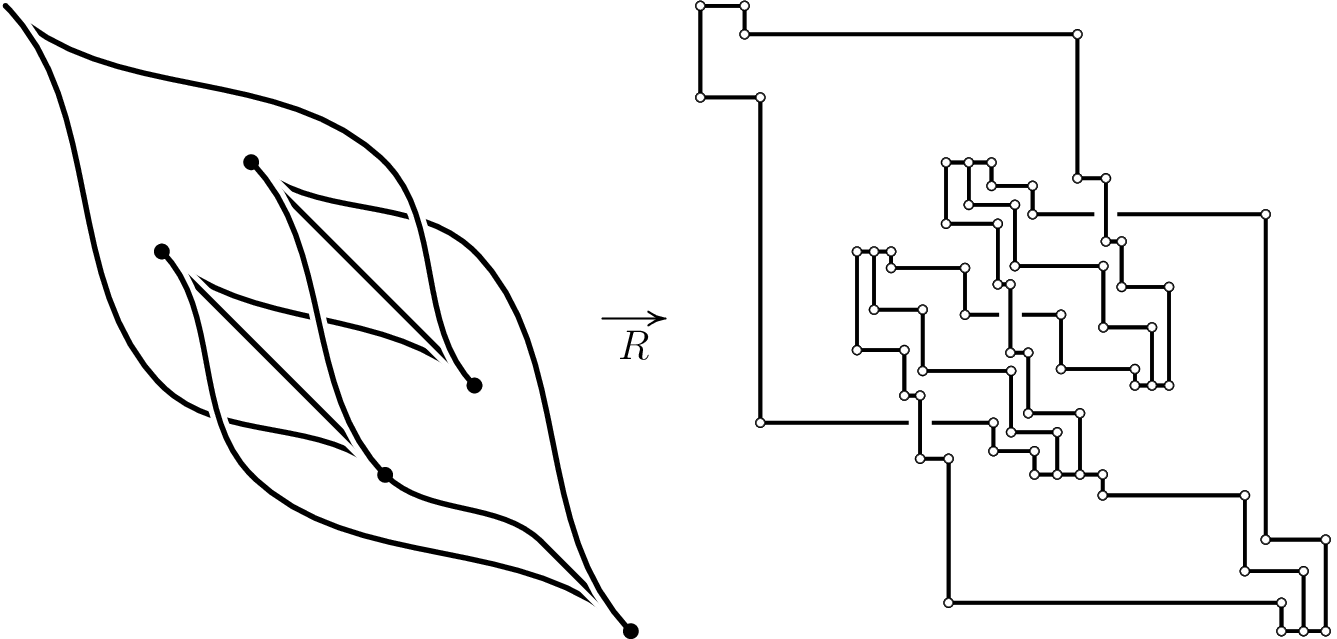}}
\caption{Approximating a front of a graph $G$ by a rectangular diagram $R(G)$}
\label{front_rect}
\end{figure}
\begin{definition}
\label{R_G}
Let $G$ be a Legendrian graph whose front is in backslash position. We say that rectangular diagram $R$ {\it approximates} $G$ if 
\begin{itemize}
\item a rectangular neighbourhood with horizontal and vertical sides is chosen for each cusp, vertex and crossing, such that the front  does  not intersect the horizontal part of its boundary; any two of these neighbourhoods does not intersect;
\item an intersection of the diagram $R$ with each neighbourhood has a prescribed form which is shown in Fig.~\ref{approx_parts}; this intersection depends only on the type of singularity of the front (a crossing, a cusp, a vertex), and in the case of a vertex the intersection depends only on the amount of the edges emerging to the left and to the right side; 
\item vertices of the diagram $R$ lying in distinct neighbourhoods do not lie on the same horizontal or vertical line;
\item a complement of the diagram $R$ to the union of these neighbourhoods is a collection of non-intersecting rectangular curves, whose corners point only to right-up or left-down direction, and  for each pair of these neighbourhoods the number of rectangular curves connecting them equals to the number of segments connecting these two neighbourhoods in the complement of the front to the union of all neighbourhoods.
\end{itemize}
Denote by $R(G)$ the set of rectangular diagrams approximating a Legendrian graph $G$. It is clear that $R(G)$ is non-empty: if we choose the neighbourhoods small enough then at least one such approximating diagram exists. An example is shown in Fig.~\ref{front_rect}.
\end{definition}

\begin{prop} Let $G$ be a graph whose front is in backslash position and $R\in R(G)$. Then the front of $G_{R}$ can be obtained from the front of $G$ by plane isotopy, blow-ups, R-moves and moves $\mathrm{II}_G$.
\end{prop}
\proof
There are two differencies between fronts of $G$ and $G_R$:
\begin{itemize}
\item $G_R$ has extra 2-valent vertices which occured from vertices of the diagram $R$;
\item all vertices of $G$ are splitted into 2- or 3-valent vertices of $G_R$.
\end{itemize}

You can add an extra 2-valent vertex easy by applying a blow-up. Then you should move it to the appropriate place using R-moves and moves $\mathrm{II}_G$ to pass cusps and crossings, see Fig.~\ref{r_move_application}.

Splitting of a vertex into 2- or 3-valent vertices can be done by applying several blow-ups.
\endproof

It is clear that any two diagrams in $R(G)$ are related by commutations and (de)stabilizations of type L. 
So the map $G\mapsto R(G)$ is the desired inverse to the map $R\mapsto G_R$ concerned in Theorem~\ref{main_result}, and to prove the theorem we only need to show that the map $G\mapsto R(G)$ is well defined and surjective. The surjectivity follows from propositions \ref{step3} and \ref{step4}. To prove the other part we are to demonstrate that any Legendrian graph after some Legendrian isotopy has a front in backslash position and that the image  of the map $G\mapsto R(G)$ is independent from the choice of the front in backslash position. This will be done in two steps: propositions \ref{step1} and \ref{step2}. 

\begin{prop}[\cite{BI}]
\label{step1}
Every Legendrian graph may be Legendrian isotoped to make its front in backslash position. Any Legendrian isotopy between two Legendrian graphs with fronts in backslash position may be continuously deformed to Legendrian isotopy between the same graphs but through Legendrian graphs with fronts in backslash position such that the deformation preserves the plane isotopy class of a front at each time. 
\end{prop}
\proof
Define the diffeomorphism $\varphi$ from the $yz$-plane to itself by $(y,z) \mapsto (y,\lambda z)$,
where $\lambda > 0$ is a positive real number, and the diffeomorphism $\psi$ as the $-\pi/4$
rotation map of the $yz$-plane. For a given Legendrian isotopy  we can choose $\lambda$ sufficiently small such that the diffeomorphism $\psi\circ\varphi$
maps all front projections during the isotopy into backslash position. Finally note that if a front $f$ is already in backslash position then a front $(\psi\circ\varphi) (f)$ is isotopic to $f$ through fronts in backslash position.
\endproof
\begin{prop}
\label{step2}
If fronts in backslash position of Legendrian graphs $G'$ and $G$ are related by some move of fronts in backslash position shown in Fig.~\ref{front_moves} and~\ref{edge_contraction} then $R(G')$ and $R(G)$ are related by elementary moves of type L.
\end{prop}
\proof
Recall that all diagrams in $R(G)$ are related by commutations and (de)stabilizations of type L. So we are done if we show that some diagram from $R(G)$ is related to some diagram in $R(G')$ by elementary moves of type L.

By Theorem \ref{front_moves_reduced} it is sufficient to consider plane isotopy of fronts, moves R, $\mathrm{II}_G$, III and a blow-up.
\begin{itemize}
\item\textit{Isotopy of fronts.} This case is obvious. Indeed, if $G'$ is obtained from $G$ by small front isotopy then by definition of $R(\cdot)$ there is a common diagram in $R(G)$ and in $R(G')$. The rest is the compactness argument.

\begin{figure}
\center{\includegraphics{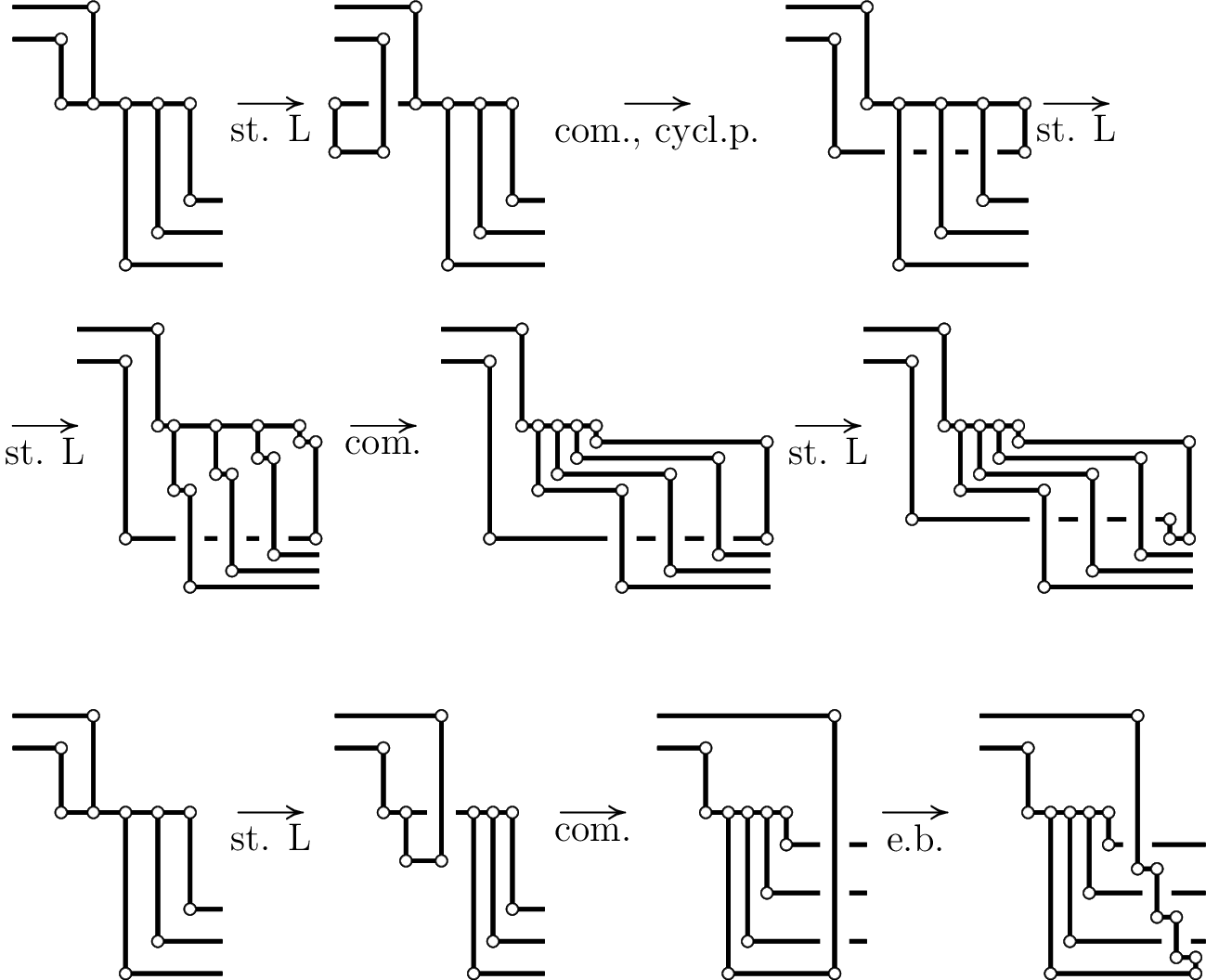}}
\caption{Translating R-move to elementary moves (of type L) of rectangular diagram}
\label{R_approx}
\end{figure}

\item\textit{R-move.} Assume that $G'$ is obtained from $G$ by R-move. We consider two cases when left-top or left-bottom branch is threw to the right side of the vertex during R-move. Other two cases differ only by $\pi$-rotation. On the picture~\ref{R_approx} it is demonstrated how to obtain some diagram belonging to $R(G')$ from some diagram belonging to $R(G)$ by elementary moves of type L.

\begin{figure}
\center{\includegraphics{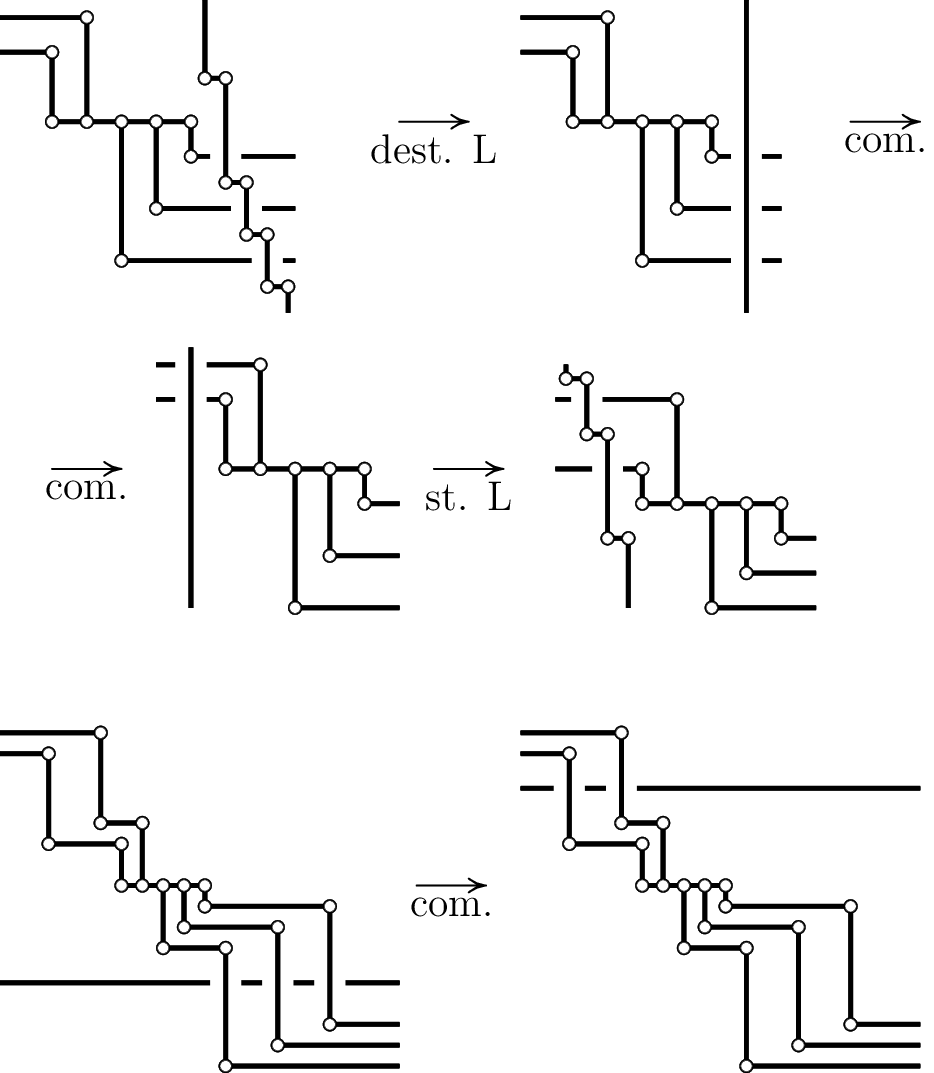}}
\caption{Translating a move $\mathrm{II}_G$ to elementary moves of type L}
\label{IIG_approx}
\end{figure}
\begin{figure}
\center{\includegraphics{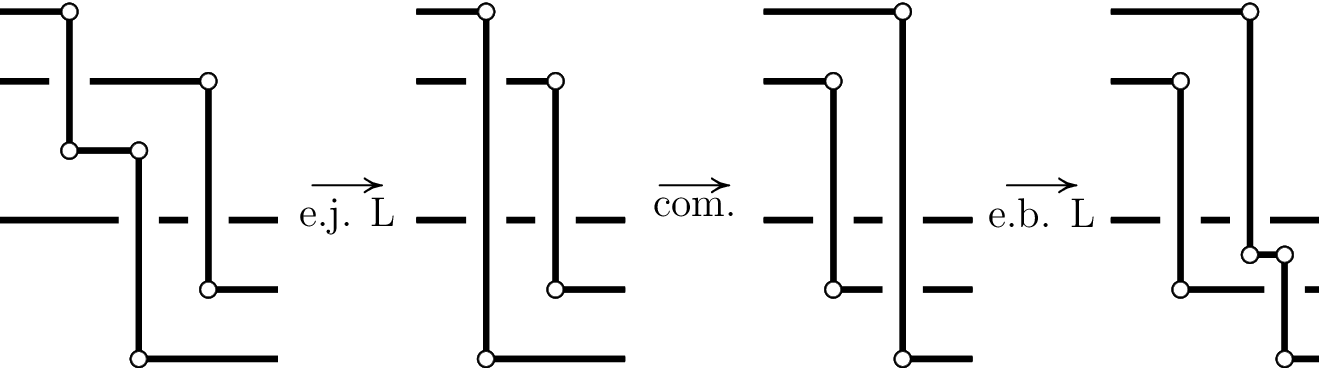}}
\caption{Translating a move III to elementary moves of type L}
\label{III_approx}
\end{figure}
\begin{figure}
\center{\includegraphics{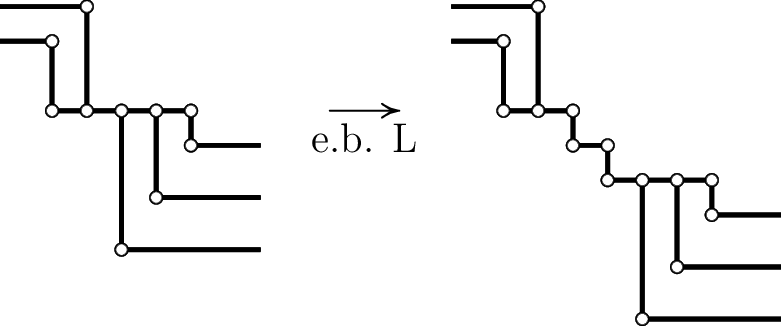}}
\caption{Translating a blow-up to elementary moves of type L}
\label{blow-up_approx}
\end{figure}

\item\textit{Move $II_G$.} We have two cases here: the line which we move through the vertex has smaller or greater slope than edges (incident to the vertex) have, see Fig.~\ref{IIG_approx}.

\item\textit{Move III.} See Fig.~\ref{III_approx}.

\item\textit{Blow-up.} In Fig.\ref{blow-up_approx} it is shown how to obtain a diagram in $R(G')$ from a diagram in $R(G)$ just by two edge breakings of type L.
\end{itemize}
\endproof

\proof[Proof of the statement that the map is well-defined.]
Suppose two Legendrian graphs $G$ and $G'$ have fronts in backslash position and are Legendrian isotopic. By proposition~\ref{step1} we can assume that they are isotopic through Legendrian graphs with fronts in backslash position. By theorem~\ref{BItheo} this isotopy splits into moves of fronts in backslash position. For each such move and for blow-up and edge contraction we proved in proposition~\ref{step2} that the image of the map $G\mapsto R(G)$ preserves. So $R(G)=R(G')$ modulo elementary moves of type L.
\endproof
\begin{prop}
\label{step3}
Suppose we have some generalized rectangular diagram $R$. It can be transformed to a diagram whose any vertical edge has two vertices by applying elementary moves of type $L$.
\end{prop}
\proof
If $R$ has a vertical edge with more than two vertices then apply an end shift of type L to top pair of vertices on this edge removing the top vertex. If some vertical edge of $R$ has only one vertex then add a vertex a little above it. After several such operations we are done.
\endproof
\begin{prop}
\label{step4}
The map $G\mapsto R(G)$ is surjective.
\end{prop}
\proof
Suppose we have some diagram $R$. We will apply type L elementary moves to $R$ to obtain a diagram of the form $R(G)$.

By the proposition \ref{step3} we can assume that every vertical edge of $R$ has exactly two vertices. Then do 4 steps of arrangements:
\begin{itemize}
\item Take some horizontal edge $h$ with at least three vertices on it. At each vertex on $h$ from right to left apply a stabilization of type L such that new small vertical edge emerges upwards from $h$. By commutations move all new small vertical edges closely to the right one. Choose a small rectangular neighbourhood (with horizontal and vertical sides) which contains all these small vertical edges such that the diagram $R$ leaves this neighbourhood from the left side. Note that such neighbourhood is illustrated in the left part of Fig.~\ref{approx_parts} in the case when there are no small vertical edges emerging downwards.
\item Suppose that a vertex $v$ lies on the same horizontal edge with exactly one another vertex, and horizontal and vertical edges emerge from $v$ to the right-bottom or the left-top. Note that in this case $v$ becomes a cusp on the front of the graph $G_R$. Then make two type L stabilizations to obtain a neghbourhood of $v$ illustrated in the center of Fig.~\ref{approx_parts}.
\item For each crossing $c$ make a type L stabilization at each end of the edge containing $c$ and (possibly) several commutations to obtain a neighbourhood of $c$ illustrated on the right part of Fig.~\ref{approx_parts}.
\item If a horizontal edge contains only one vertex then apply a type L stabilization at it. Choose a rectangular neighbourhood which contains the small vertical edge and intersects the small horizontal edge.
\end{itemize}

Note that in these neighbourhoods the diagram has a prescribed form shown in Fig.~\ref{approx_parts}, and outside neighbourhoods: the diagram has no crossings, the diagram has only edges with two vertices, and all angles are left-bottom and right-top. 
\endproof

\section{Fence diagrams}

\begin{figure}
\center{\includegraphics{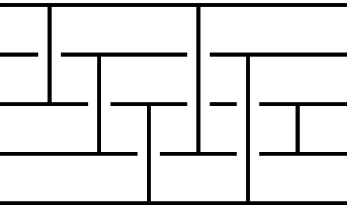}}
\caption{Fence diagram}
\label{fence_diagram}
\end{figure}

\begin{definition}[\cite{Rud}]
A {\it fence diagram} consists of horizontal segments (on the plane) called {\it posts} and vertical segments called {\it wires} such that:
\begin{itemize}
\item all posts are related by $y$-shift;
\item ends of wires are the interior points of posts;
\item all wires lie on distinct horizontal levels.
\end{itemize}
Two fence diagrams are assumed to be the same if they are combinatorially equivalent (as in the definition \ref{GRD}).
At the crossings we always assume that the wire is  above the post. An example is shown in Figure \ref{fence_diagram}.
\end{definition}

\begin{figure}
\center{\includegraphics{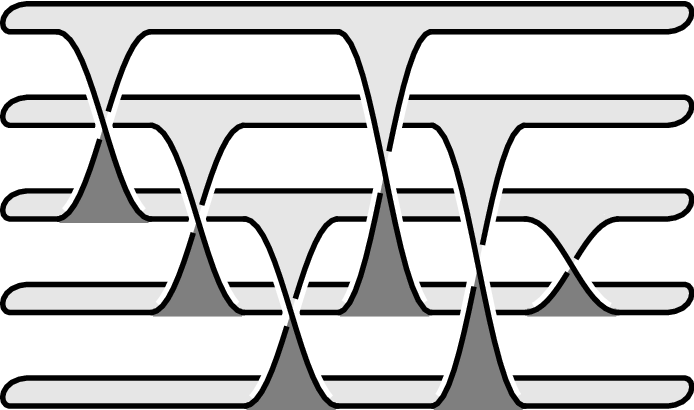}}
\caption{Quasipositive surface}
\label{quasipositive_surface}
\end{figure}

From a fence diagram Rudolph in \cite{Rud} constructs a surface (with boundary) embedded in $\mathbb{R}^3$: each post corresponds to a horizontal disc and each wire --- to a positive band joining two discs. The surface obtained by this construction is called {\it quasipositive}. An example is shown in Figure \ref{quasipositive_surface}.

\begin{figure}
\center{\includegraphics{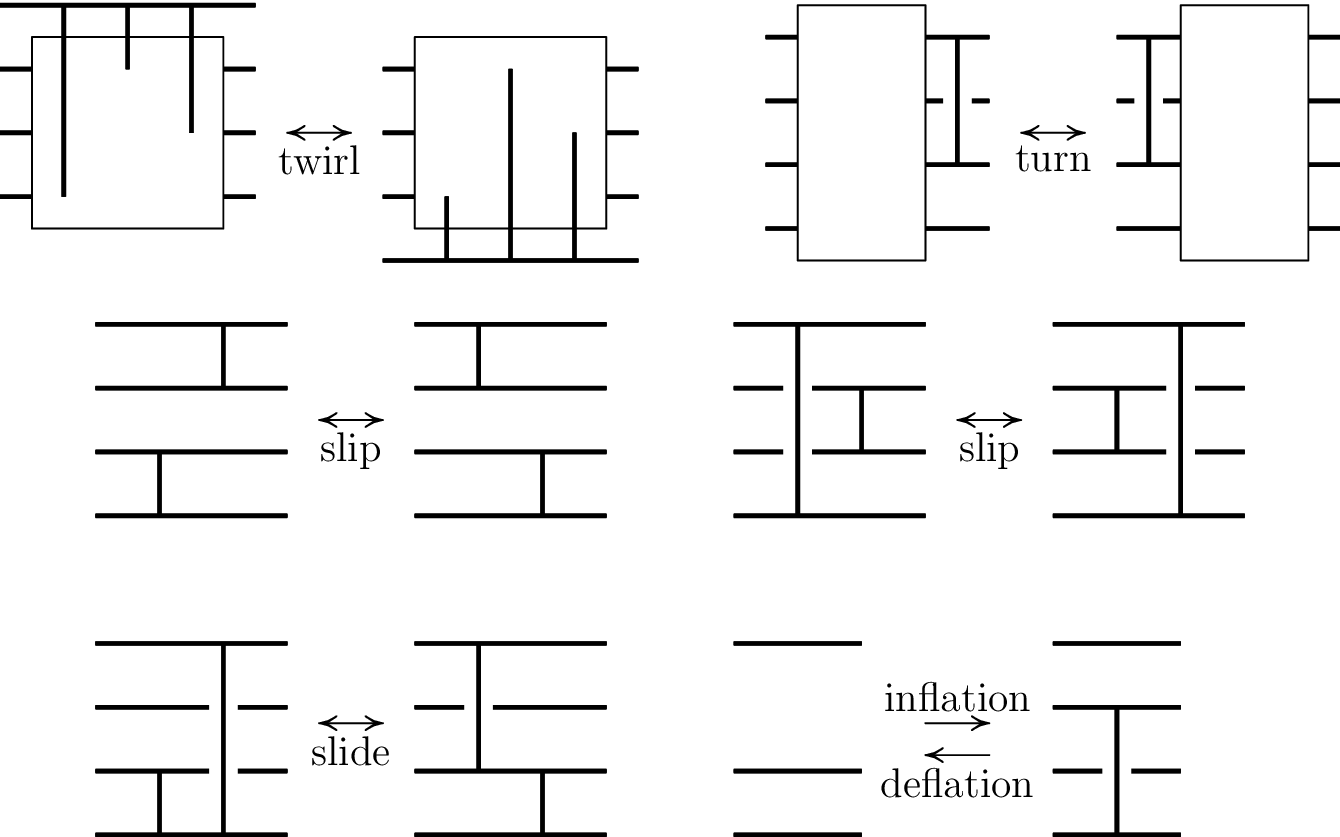}}
\caption{Fence moves}
\label{fence_moves}
\end{figure}

To classify isotopy classes of quasipositive surfaces Rudolph in \cite{Rud2} introduced {\it fence moves} (see Fig. \ref{fence_moves}):
\begin{itemize}
\item A {\it twirl} consists in moving top or bottom post to the opposite side taking along all wires connecting this post to others;
\item A {\it turn} consists in moving the leftmost or rightmost wire to the opposite side;
\item A {\it slip} consists in exchange of horizontal positions of two neighboring wires provided that they do not interleave;
\item Suppose we have two neighboring wires such that the left one connects posts on vertical levels $a$ and $b$ and the right one connects posts on levels $a$ and $c$ where $a<b<c<a$ in the cyclic order. A {\it slide} consists in substituting for these wires such two wires on the same horizontal levels that the left one connects posts on vertical levels $b$ and $c$ and the right one connects posts on levels $b$ and $a$.
\item An {\it inflation} consists in adding one post on any vertical level and connecting it with any another post by one wire. The inverse operation is called {\it deflation}.
\end{itemize}

Rudolph asked in \cite{Rud2}, is it true that two fence diagrams are related by fence moves if and only if the corresponding quasipositive surfaces are isotopic? Baader and Ishikawa in \cite{BI} answered in the negative. We briefly discuss their approach. 

They construct a map from 3-valent Legendrian graphs to fence diagrams modulo fence moves. This is done as in the construction of the map $R(G)$ (see definition \ref{R_G}): they approximate a front in backslash position (see definition \ref{backslash_position}) by a fence diagram and show that if the front is changed by a plane isotopy or a front move then the approximating fence diagram can be changed by fence moves.

They also notice that if a fence diagram corresponds to a Legendrian link then fence moves preserve the Legendrian type of that link. In the final they provide two isotopic quasipositive surfaces whose fence diagrams have distinct corresponding Legendrian links (distinguished by a rotation number).

\begin{remark}
Also in \cite{BI} authors noticed that a Legendrian ribbon (a surface defined in \ref{ribbon}) is quasipositive. For other connections between quasipositive surfaces and contact geometry see, for example, \cite{Rud3}.
\end{remark}

Our aim in this section is to factorize the map (constructed by Baader and Ishikawa) to Legendrian graphs modulo blow-up and to prove that this new map is a bijection. For factorizing the map it is sufficient to show a little: that if a front is changed by a blow-up then the approximating fence diagram can be changed by fence moves. But instead we will construct a bijective map from fence diagrams modulo fence moves to generalized rectangular diagrams modulo elementary moves of type L.

\begin{theorem}
\label{fence_main_result}
Denote by $3\text{-}LG$ the set of Legendrian graphs modulo Legendrian isotopy whose vertices have valence $2$ or $3$, by $FD$ --- fence diagrams modulo fence moves, by $LR$ ("Legendrian Ribbons") --- Legendrian graphs modulo Legendrian isotopy and blow-up, by $GRD_L$ --- generalized rectangular diagrams modulo elementary moves of type L. 

Let $3\text{-}LG\to FD$ be the map constructed in \cite{BI}, $3\text{-}LG\to LR$ --- a natural map, $GRD_L\to LR$ --- the map $R\mapsto G_R$ defined in \ref{mapG}. Then there exists a bijective map $FD\to GRD_L$ such that the following diagram is commutative:
$$\begin{matrix}
3\text{-}LG & \rightarrow & FD \\
\downarrow & & \downarrow \\
LR & \leftarrow & GRD_L
\end{matrix}$$
\end{theorem}

Before proving the theorem we introduce two definitions.

\begin{figure}
\center{\includegraphics{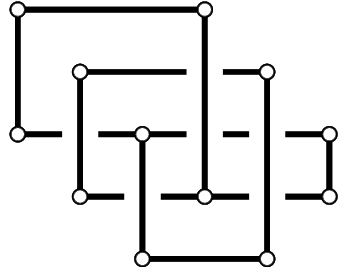}}
\caption{A rectangular diagram obtained from the fence diagram in Figure \ref{fence_diagram}}
\label{R_F}
\end{figure}

\begin{definition}
Denote by $R(F)$ the generalized rectangular diagram whose vertices are the ends of wires of the fence diagram $F$. An example is shown in Figure \ref{R_F}.
\end{definition}

\begin{definition}
Let $R$ be a generalized rectangular diagram such that each vertical edge contains two vertices. Denote by $F(R)$ a fence diagram whose wires are the vertical edges of $R$ and whose posts contain horizontal edges of $R$.

If $R$ has a vertical edge with more than two vertices then apply an end shift of type L to top pair of vertices on this edge removing the top vertex. If some vertical edge of $R$ has only one vertex then add a vertex a little above it. After several such operations we obtain a diagram with vertical edges having only two vertices. Denote the fence diagram corresponding to the obtained rectangular diagram also by $F(R)$.
\end{definition}

\noindent{\it Proof of Theorem \ref{fence_main_result}}. Note that $F(R(F))=F$. We will prove that $F\mapsto R(F)$ is the desired bijection. It follows from two lemmas.

\begin{lemma}
If fence diagram $F'$ is obtained from $F$ by a fence move then $R(F')$ can be obtained from $R(F)$ by elementary moves of type L.
\end{lemma}
\proof We consider in each case of a fence move which elementary moves are to apply:
\begin{itemize}
\item twirl: a cyclic permutation of horizontal edges;
\item turn: a cyclic permutation of vertical edges;
\item slip: a commutation of vertical edges;
\item slide: a combination of an end shift and the inverse of the end shift;
\item inflation: two vertex additions.
\end{itemize}
\endproof
\begin{lemma}
If generalized rectangular diagram $R'$ is obtained from $R$ by elementary move of type L then $F(R')$ can be obtained from $F(R)$ by fence moves.
\end{lemma}
\proof By Theorem \ref{elementary_basis} it is sufficient to consider cases of commutations, an end shift of type L and a vertex addition. In each case we provide which fence moves are to apply:
\begin{itemize}
\item {\it A commutaion of vertical edges}. Several slips.

\begin{figure}
\center{\includegraphics{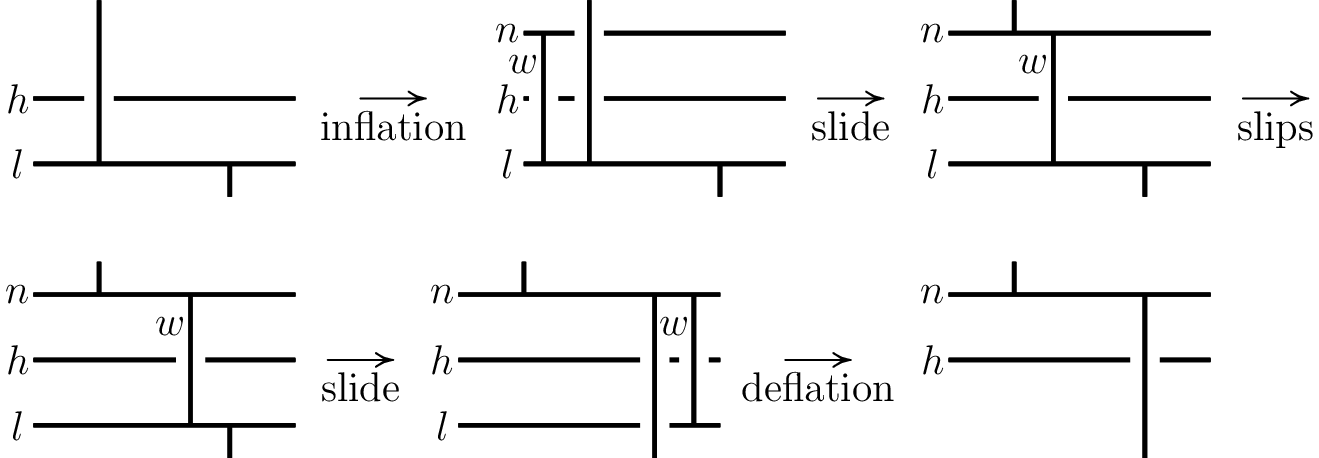}}
\caption{Translating horizontal commutation to fence moves}
\label{horizontal_commutation_fence}
\end{figure}

\item {\it A commutation of horizontal edges}. Denote a post corresponding to the higher commuting edge of $R$ by $h$ and the lower post by $l$. Make several turns to make the wires attached to the post $l$ being on the left to the wires attached to the post $h$. Make an inflation: a new post $n$ will be a little higher than the post $h$ and a new wire $w$ will be attached to the post $l$ on the left to wires attached to this post. Apply a slide to $w$ and the neighbouring wire on $l$. As a result the wire $w$ moves on the right and the other wire becomes attached to the post $n$. Make slips to move $w$ to the next wire on $l$. Continue in the same way. In the end we obtain that all wires initially attached to $l$ become attached to $n$ and the wire $w$ is the single wire attached to $l$. So make a deflation to remove $w$ and $l$. A horizontal commutation is almost done: make several turns to move wires to initial horizontal positions. See Fig. \ref{horizontal_commutation_fence}.

\item {\it Vertical end shift}. Some number (may be zero) of slides and, may be, an inflation.

\begin{figure}
\center{\includegraphics{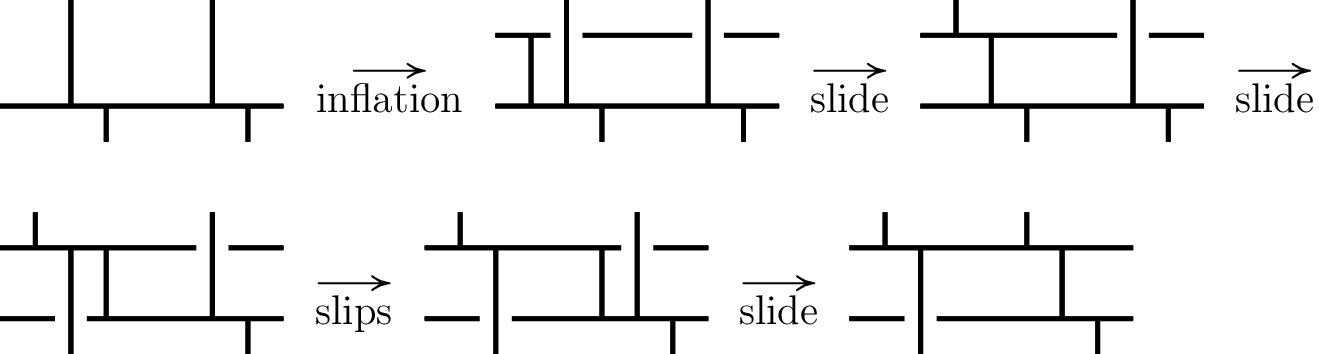}}
\caption{Translating horizontal end shift to fence moves}
\label{horizontal_end_shift_fence}
\end{figure}

\item {\it Horizontal end shift}. Consider the case of horizontal end shift where the left vertex is removed. The other case can be considered by applying rotation by $\pi$. By construction of $F(R)$ near the vertex to remove we have one or two wires. Make an inflation: add a post a little higher and add a wire on the left near the vertex to remove. Then apply, respectively, one or two slides. Then apply several slips to move new wire to the other vertex concerned in the end shift. There we also have one or two wires. Apply, respectively, zero or one slide and we are done. An example is shown in Figure \ref{horizontal_end_shift_fence}.

\item {\it A vertex addition on the vertical edge.} An inflation and, in some cases, one slide (in the other cases nothing else to do).

\item {\it A vertex addition on the horizontal edge.} An inflation.
\end{itemize}
\endproof
\endproof

\begin{corollary}
\label{fence_leg}
A map $3\text{-}LG\to FD$ constructed in \cite{BI} induces a bijection $LR\to FD$. So Legendrian graphs modulo Legendrian isotopy and blow-up are in one-to-one correspondence with fence diagrams modulo fence moves.
\end{corollary}

\section{Spatial graphs}

\begin{figure}
\center{\includegraphics{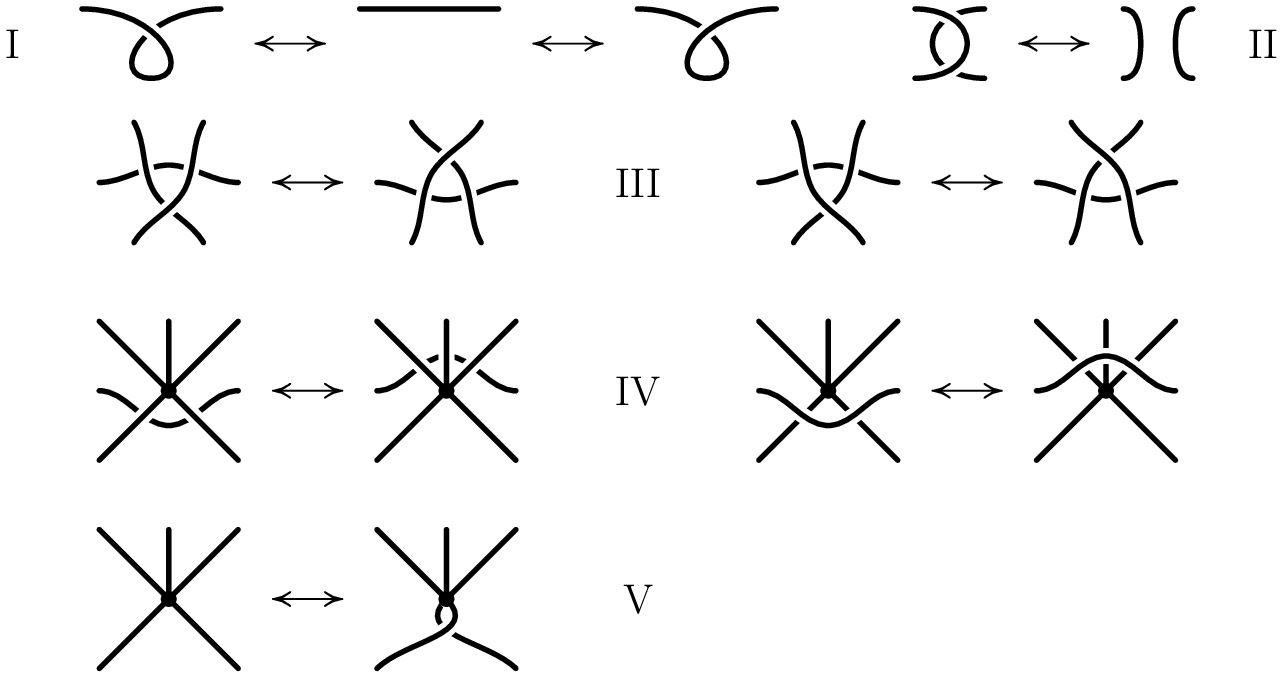}}
\caption{Reidemeister moves}
\label{Reidemeister_moves}
\end{figure}

In \cite{Kauffman} the analogue of Reidemeister theorem for spatial graphs was proved: that two spatial graphs are isotopic if and only if their generic plane projections are related by plane isotopy and Reidemeister moves shown in Figure \ref{Reidemeister_moves}.

Using this theorem one can prove a variant of Theorem \ref{main_result} for spatial graphs:

\begin{prop}
Consider any generalized rectangular diagram as a plane diagram of some spatial graph. This correspondence induces a bijection between spatial graphs modulo isotopy and blow-up and generalized rectangular diagrams modulo elementary moves.
\end{prop}
\noindent{\it Sketch of the proof.} In the proof of Theorem \ref{main_result} we already proved that elementary moves of type L may cause only an isotopy, a blow-up or an edge contraction of the corresponding graph. Apply to all diagrams a horizontal symmetry: then moves of type L become moves of type N and the topological type of the graph changes as if only the orientarion of the ambient space is reversed. So elementary moves of type N also do not alter the class of a spatial graph. And we get a map from rectangular diagrams modulo elementary moves to spatial graphs modulo isotopy and blow-up.

To prove that this map is a bijection consider the same argument as in Theorem \ref{main_result}. Suppose we have a plane diagram of a spatial graph. Rotate all its crossings so that the overpass has vertical tangent and the underpass has horizontal tangent at each crossing. Then approximate the obtained plane diagram by generalized rectangular diagram. We are done after we show that:
\begin{itemize}
\item A class of the obtained rectangular diagram (modulo elementary moves) is independent of the way of rotating crossings.
\item If two plane diagrams differ by a Reidemeister move then the corresponding rectangular diagrams are related by elementary moves.
\end{itemize}
Proofs are easy and do not involve new ideas. 
\endproof

{\sc
Laboratory of Quantum Topology

Chelyabinsk State University

Brat'ev Kashirinykh street 129 

Chelyabinsk 454001, Russia.

\smallskip
\smallskip
\smallskip

Dept. of Mechanics \& Mathematics

Moscow State University

Leninskije gory street 1

Moscow 119234, Russia}

\smallskip
\smallskip
\smallskip

\texttt{0x00002a@gmail.com}

\end{document}